\documentclass[11pt]{amsart}

\usepackage{geometry}
\usepackage{todonotes}
\usepackage{amsfonts, amssymb, amscd}
\numberwithin{equation}{section}

\usepackage[symbol]{footmisc}

\usepackage{bm}
\usepackage{verbatim}
\usepackage{mathrsfs}
\usepackage{graphicx}
\usepackage{tikz-cd}
\usepackage{subcaption}
\usepackage{listings}
\usepackage{subfiles}
\usepackage[toc,page]{appendix}
\usepackage{mathtools}
\usepackage{comment}
\usepackage{enumerate}
\usepackage{enumitem}
\usepackage[all]{xy}

\usepackage{graphicx}
\graphicspath{{images/}}

\usepackage{appendix}
\usepackage{hyperref}
\hypersetup{
    colorlinks=true,
    citecolor=red,
    linkcolor=blue,
    filecolor=magenta,      
    urlcolor=red,
}
\newcommand{\bb}{\bm{b}}
\newcommand{\Mm}{{\bf{M}}}

\newcommand{\NN}{{\bf{N}}}
\newcommand{\Dd}{{\bf{D}}}

\newcommand{\Qq}{\mathbb{Q}}

\newcommand{\Rr}{\mathbb{R}}

\newcommand{\Center}{\operatorname{center}}

\newcommand{\Exc}{\operatorname{Exc}}

\newcommand{\glct}{\operatorname{glct}}

\newcommand{\rk}{\operatorname{rank}}

\newcommand{\Weil}{\operatorname{Weil}}
\newcommand{\lct}{\operatorname{lct}}

\newcommand{\Supp}{\operatorname{Supp}}
\newcommand{\Ngklt}{\operatorname{Ngklt}}

\newcommand{\mult}{\operatorname{mult}}

\newcommand{\Aa}{{\bf{A}}}

\newcommand{\Pic}{\mathrm{Pic}}

\newtheorem{thm}{Theorem}[section]

\newtheorem{cor}[thm]{Corollary}
\newtheorem{lem}[thm]{Lemma}

\newtheorem{claim}[thm]{Claim}

\theoremstyle{definition}
\newtheorem{defn}[thm]{Definition}
\newtheorem{ques}[thm]{Question}
\theoremstyle{definition}

\newtheorem{deflem}[thm]{Definition-Lemma}
\newtheorem{ex}[thm]{Example}

\newtheorem{cond}[thm]{Condition}

\theoremstyle{definition}

\begin{document}

\title{On generalized lc pairs with $\bb$-log abundant nef part}
\author{Junpeng Jiao, Jihao Liu, and Lingyao Xie}

\address{Department of Mathematics, The University of Utah, Salt Lake City, UT 84112, USA}
\email{jiao@math.utah.edu}

\address{Department of Mathematics, Northwestern University, 2033 Sheridan Rd, Evanston, IL 60208}
\email{jliu@northwestern.edu}

\address{Department of Mathematics, The University of Utah, Salt Lake City, UT 84112, USA}
\email{lingyao@math.utah.edu}

\subjclass[2020]{14E30,14C20.14E05,14J17}
\date{\today}

\begin{abstract}
We study the behavior of generalized lc pairs with $\bb$-log abundant nef part, a meticulously designed structure on algebraic varieties. We show that this structure is preserved under the canonical bundle formula and sub-adjunction formulas, and is also compatible with the non-vanishing conjecture and the abundance conjecture in the classical minimal model program. 
\end{abstract}

\maketitle

\tableofcontents

\section{Introduction}

We work over the field of complex numbers $\mathbb C$.

The classical birational geometry and the minimal model program originated with the study of smooth complex projective varieties $X$ and their canonical divisors $K_X$. It turns out that in order to understand smooth varieties, it is also necessary to study varieties with mild singularities and several natural structures on these varieties. More precisely, starting from a smooth projective variety, the minimal model program predicts that after finitely many flips and divisorial contractions, we will reach a birational model $X'$ of $X$ with terminal singularities, such that
\begin{enumerate}
    \item either $X'$ is a good minimal model (i.e. $K_{X'}$ is semi-ample) and the induced morphism $f:X'\to Z:={\mathbf{Proj}}(R(K_X))$ is the Iitaka fibration, or
    \item there exists a morphism $f:X'\to Z$ which is a Mori fiber space, and in particular, $\dim Z<\dim X'$ and $-K_{X'}$ is ample$/Z$.
\end{enumerate}
Thus the study of these two kinds of fibrations is of fundamental importance in the minimal model program.
In the study of these fibrations  $X\rightarrow Z$, a structure called the \emph{canonical bundle formula} naturally appears. This was first studied for minimal elliptic fibrations  in \cite{Kod63}. More precisely,we have the following.  In case (1), we can always write
$$K_X\sim_{\Qq}f^*(K_Z+B_Z+\Mm^Z_Z)$$
for some $\Qq$-divisor $B_Z$ measuring the singularities of the fibration (defined by certain lc thresholds) and a nef $\bb$-divisor $\Mm^Z$ measuring the variation of the morphism.
In case (2), for any general ample$/Z$ $\Qq$-divisor $0\leq B\sim_{\Qq}-K_X$, one can write
$$K_X+B\sim_{\Qq}f^*(K_Z+B_Z+\Mm^Z_Z)$$
for the pair $(X,B)$ and thus we reduce to a situation similar to case (1). The above formulas are called \emph{canonical bundle formulas}. In fact, it is known that for any pair $(X,B)$ that is lc over the generic point of $Z$ and $K_X+B\sim_{\Rr,Z}0$, we also have canonical bundle formula  inducing a similar structure $(Z,B_Z,\Mm^Z)$ (cf. \cite{Kaw98,Amb99,PS09}).

It is interesting to ask whether we can understand the structure of $(Z,B_Z,\Mm^Z)$ more precisely. Prokhorov-Shokurov conjectured that $\Mm^Z$ is a $\bb$-semi-ample \cite[Conjecture 7.13]{PS09}, however this is known to be a very difficult conjecture and is not even known in the case when $\dim X=3$. A weaker conjecture is that $\Mm^Z$ is $\bb$-abundant. This was proven by Ambro \cite[Theorem 3.3]{Amb05} when $(X,B)$ is klt, and by Fujino-Gongyo \cite[Theorem 1.1]{FG14} when $(X,B)$ is lc.

In 2014, Birkar and Zhang \cite{BZ16} introduced a more general structure for algebraic varieties, called \emph{generalized pair} (\emph{g-pair} for short):
\begin{defn}[{cf. Definition \ref{defn: g-pairs}}]\label{defn: g-pair intro}
A generalized pair (g-pair for short) $(X,B,\Mm)/U$ consists of a projective morphism $X\rightarrow U$ from a normal quasi-projective variety to a variety, an $\Rr$-divisor $B\ge 0$ on $X$, and an NQC$/U$ $\bb$-divisor $\Mm$ over $X$, such that $K_X+B+\Mm_X$ is $\Rr$-Cartier.
\end{defn}
G-pairs play an important role in many recent developments in birational geometry, such as the effective Iitaka fibration \cite{BZ16}, the proof of the BAB conjecture \cite{Bir19,Bir21a}, and the connectedness principles \cite{HH19,Bir20b,FS20}. The classical minimal model program also works well for all glc g-pairs \cite{BZ16,HL21}. We refer the reader to \cite{Bir20a} for a survey on other recent progress.

\medskip

\noindent\textbf{Canonical bundle formula for generalized pairs}. It is immediate from the definition that any structure $(Z,B_Z,\Mm^Z)$ deduced from the canonical bundle formula is a g-pair. Therefore, the canonical bundle formula allows us to apply the structure of g-pairs in lower dimensions to study the behavior of pairs in higher dimensions. Following this philosophy, one question naturally arises for inductive purposes: can we find a canonical bundle formula between g-pairs? Fortunately for us, we already have a partial positive answer to this question:

\begin{thm}[{cf. \cite[Theorem 1.2]{HL19}; see also \cite[Theorem 2.20]{FS20}, \cite[Theorem 1.4]{Fil20}}]\label{thm: gcbf}
Let $(X,B,\Mm)/U$ be a glc (resp. gklt) g-pair and $f: X\rightarrow Z$ a projective surjective morphism$/U$, such that $K_X+B+\Mm_X\sim_{\Rr,Z}0$. Then there exists a glc (resp. gklt) g-pair $(Z,B_Z,\Mm^Z)/U$ such that
$$K_X+B+\Mm_X\sim_{\Rr}f^*(K_Z+B_Z+\Mm^Z_Z).$$
Moreover, if $(X,B,\Mm)$ is a $\Qq$-g-pair, then we may choose $(Z,B_Z,\Mm^Z)/U$ to be the g-pair induced by a canonical bundle formula\footnote{A canonical bundle formula for a projective surjective morphism $f$ is given as the composition of a Fujino-Gongyo type canonical bundle formula \cite{FG12} and a Kodaira type canonical bundle formula \cite{Kod63} via the Stein factorization of $f$. We refer the reader to Subsection \ref{subsection: scbf} formal definitions.}  of $f: (X,B,\Mm)\rightarrow Z$.
\end{thm}

A massive disadvantage for Theorem \ref{thm: gcbf} is the following: unlike the usual canonical bundle formulas for varieties and pairs, we cannot expect either the semi-ampleness or the $\bb$-abundance property of $\Mm^Z$ in general. Therefore, we ask the following questions:
\begin{ques}\label{ques: special gcbf}
Let $(X,B,\Mm)/U$ and $(Z,B_Z,\Mm^Z)/U$ be g-pairs as in Theorem \ref{thm: gcbf}. Under what additional conditions on $\Mm$, will
\begin{enumerate}
    \item $\Mm^Z$ be $\bb$-semi-ample$/U$?
    \item $\Mm^Z$ be $\bb$-abundant$/U$?
\end{enumerate}
Naturally, we will insist that these additional conditions are always satisfied by the case $\Mm=\bf{0}$ (so that the usual pair case is included).
\end{ques}
Of course, $\Mm$ is necessarily semi-ample for (1) and necessarily abundant for (2), otherwise we get a contradiction by simply considering the identity morphism $f: X\rightarrow X:=Z$.  In this case, it is proven by Filipazzi (\cite[Chapter 6, Theorem 10]{Fil19}; see also \cite[Theorem 7.3]{Fil20}) that if we have a positive answer for (1) when $\Mm=\bf{0}$, then we have a positive answer for (1) for all semi-ample$/U$ nef part $\Mm$, provided that $(X,B,\Mm)$ is a $\Qq$-g-pair. 

In this paper, we provide a satisfactory answer for (2):

\begin{thm}\label{thm: log abundance gcbf}
Let $(X,B,\Mm)/U$ be a glc g-pair and $f: X\rightarrow Z$ a projective surjective morphism$/U$, such that $K_X+B+\Mm_X\sim_{\Rr,Z}0$ and $\Mm$ is $\bb$-log abundant$/U$ (see Definition \ref{defn: b log abundance}) with respect to $(X,B,\Mm)$. Then we have a canonical bundle formula
$$K_X+B+\Mm_X\sim_{\Rr}f^*(K_Z+B_Z+\Mm^Z_Z)$$
such that $\Mm^Z$ is NQC$/U$ and $\bb$-log abundant$/U$ with respect to $(Z,B_Z,\Mm^Z)$.
\end{thm}
We note that Ambro \cite[Theorem 3.3]{Amb05} proves Theorem \ref{thm: log abundance gcbf} when $\Mm=\bf{0}$, $f$ is a contraction, and $(X,B)$ is a klt $\Qq$-pair (see \cite[Theorem 1.2]{Hu20} for a more general case). We also note that the assumption ``$\bb$-log abundant" in Theorem \ref{thm: log abundance gcbf} is a natural assumption in the study of lc and glc singularities, and is equivalent to ``$\bb$-abundant" for klt and gklt singularities. We refer the reader to Definition \ref{defn: b abundance} for more details. 

\medskip

\noindent\textbf{Sub-adjunction formula for generalized pairs}. Kawamata's sub-adjunction formula states that if $W$ is an lc center of an lc $\Qq$-pair $(X,B)$, then
$$K_W+B_W+\Mm^W_W\sim_{\Qq}(K_X+B)|_W$$
for some $\Qq$-glc g-pair $(W,B_W,\Mm^W)$ (cf. \cite{Kaw98}). This formula provides important information in the study of high codimensional lc centers. Thanks to work of J. Han and W. Liu, we also have a sub-adjunction formula for generalized pairs. That is, given a glc g-pair $(X,B,\Mm)/U$ and the normalization of a glc center $W$ of $(X,B,\Mm)$, we have
\begin{equation}\label{equ: sub-adjunction gpair}
    K_W+B_W+\Mm^W_W\sim_{\Rr}(K_X+B+\Mm_X)|_W
\end{equation}
for some glc g-pair $(W,B_W,\Mm^W)/U$ (\cite[Theorem 5.1]{HL19}). Despite its technical appearance, (\ref{equ: sub-adjunction gpair}) is very useful for inductive purposes and plays a crucial role in the proof of the cone theorem of glc g-pairs \cite[Proof of Lemma 5.14]{HL21}.  An immediate application of Theorem \ref{thm: log abundance gcbf} is the sub-adjunction formula for glc g-pairs with $\bb$-log abundant nef part:

\begin{thm}\label{thm: sub-adjunction log abundant}
Let $(X,B,\Mm)/U$ be a glc g-pair such that $\Mm$ is $\bb$-log abundant$/U$ with respect to $(X,B,\Mm)$. Let $W$ be the normalization of a glc center of $(X,B,\Mm)$ of dimension $\geq 1$. Then there exists a glc g-pair $(W,B_W,\Mm)/U$ given by the sub-adjunction
$$K_W+B_W+\Mm^W_W\sim_{\Rr}(K_X+B+\Mm_X)|_W$$
such that $\Mm^W$ is $\bb$-log abundant$/U$ with respect to $(W,B_W,\Mm^W)$.
\end{thm}

It is important to notice that Theorems \ref{thm: log abundance gcbf} and \ref{thm: sub-adjunction log abundant} will allow us to apply induction on dimension for glc g-pairs with $\bb$-log abundant nef part.

\begin{comment}
\begin{cor}\label{cor: gklt abundance gcbf}
Let $(X,B,\Mm)/U$ be a gklt g-pair and $f: X\rightarrow Z$ a projective surjective morphism$/U$, such that $K_X+B+\Mm_X\sim_{\Rr,Z}0$ and $\Mm$ is $\bb$-abundant$/U$ with respect to $(X,B,\Mm)$. Then we have a canonical bundle formula
$$K_X+B+\Mm_X\sim_{\Rr}f^*(K_Z+B_Z+\Mm^Z_Z)$$
for some glc g-pair, such that $\Mm^Z$ is $\bb$-abundant$/U$ with respect to $(Z,B_Z,\Mm^Z)$. 

Moreover, if $(X,B,\Mm)$ is a $\Qq$-g-pair and $f$ is a contraction, then we may choose $(Z,B_Z,\Mm^Z)/U$ to be the g-pair induced by a canonical bundle formula of $f: (X,B,\Mm)\rightarrow Z$.
\end{cor}
\end{comment}

\medskip

\noindent\textbf{Non-vanishing and log abundance for generalized pairs}. Many standard conjectures in birational geometry for pairs have analogues for g-pairs. Some of these analogues hold (e.g. the effective birationality \cite{BZ16}, ACC for glc thresholds and global ACC \cite{BZ16}, the boundedness of complements (\cite{Bir19,Che20}), the connectedness principle \cite{Bir20b,FS20}, DCC for volumes \cite{Bir21b}, and the cone theorem and the existence of flips \cite{HL21}). 

However, not all analogues hold in a satisfactory fashion. The non-vanishing conjecture and the abundance conjecture both fail for generalized pairs even in dimension $1$ by considering an elliptic curve and a numerically trivial non-torsion divisor on it. In fact, even if we consider the numerical non-vanishing conjecture (\cite[Conjecture 1.2]{HL20}; see also \cite[Section 1]{LP20a}) and the numerical abundance conjecture (\cite[Section 1]{LP20a}) for generalized pairs, there are counter-examples in general (cf. \cite[Examples 6.1, 6.2]{LP20a}, \cite[Examples 1.3, 3.14]{HL20}). Nevertheless, if we add some conditions to the nef part $\Mm$, these analogical conjectures may hold. As an application of Theorem \ref{thm: gcbf}, we prove the non-vanishing and log abundance for glc g-pairs with $\bb$-log abundant nef part assuming the corresponding conjectures for klt pairs:

\begin{thm}\label{thm: klt abundant to glc good nef abundant}
Let $d$ be a positive integer, and $(X,B,\Mm)/U$ a glc g-pair of dimension $d$, such that $\Mm$ is $\bb$-log abundant$/U$ with respect to $(X,B,\Mm)$. Assume the non-vanishing conjecture for all klt pairs in dimension $\leq d$. Then
\begin{enumerate}
    \item if $K_X+B+\Mm_X$ is pseudo-effective$/U$, then $|K_X+B+\Mm_X/U|_{\Rr}\not=\emptyset$, and
    \item if we further assume the abundance conjecture for all klt pairs in dimension $\leq d$, then $K_X+B+\Mm_X$ is log abundant$/U$ with respect to $(X,B,\Mm)$. In particular, $\kappa_{\iota}(X/U,K_X+B+\Mm_X)=\kappa_{\sigma}(X/U,K_X+B+\Mm_X)$.
\end{enumerate}
\end{thm}
 
Since the non-vanishing and abundance hold for klt pairs in dimension $3$, we immediately have the following corollary:
 \begin{cor}\label{cor: log abundant m dim 3 abundance}
Let $(X,B,\Mm)/U$ be a glc g-pair of dimension $\leq 3$, such that $\Mm$ is $\bb$-log abundant$/U$ with respect to $(X,B,\Mm)$. Then $K_X+B+\Mm_X$ is $\bb$-log abundant$/U$ with respect to $(X,B,\Mm)$. In particular, $\kappa_{\iota}(X/U,K_X+B+\Mm_X)=\kappa_{\sigma}(X/U,K_X+B+\Mm_X)$.
\end{cor}
Theorems \ref{thm: log abundance gcbf} and \ref{thm: klt abundant to glc good nef abundant} indicate that the category of g-pairs with $\bb$-log abundant nef part is a meticulous structure which epitomizes the advantages from both the nice numerical property of usual pairs and the induction convenience of g-pairs. Therefore, we expect this category of g-pairs to play an important role in the future studies of both pairs and generalized pairs.

The following example is complementary to Theorems \ref{thm: log abundance gcbf} and \ref{thm: klt abundant to glc good nef abundant}, which shows that we cannot replace the assumption ``$\bb$-log abundant" with ``$\bb$-abundant":
\begin{ex}[{=Example \ref{ex: counter-example not b-log abundant}}]\label{ex: intro example ab not log ab}
There exists a contraction $f: X\rightarrow E$ from a smooth projective surface to a smooth elliptic curve, and a projective glc $\Qq$-g-pair $(X,B,\Mm)$, such that
\begin{enumerate}
    \item $\Mm=\overline{\Mm_X}$ is $\bb$-abundant but is not $\bb$-log abundant with respect to $(X,B,\Mm)$,
    \item $K_X+B+\Mm_X$ is pseudo-effective but not effective. In particular, $K_X+B+\Mm_X$ is not abundant, and
    \item $K_X+B+\Mm_X\sim_{\Qq,E}0$ but $\Mm^E$ is not $\bb$-abundant for any canonical bundle formula
    $$K_X+B+\Mm_X\sim_{\Qq}f^*(K_E+B_E+\Mm^E_E).$$
\end{enumerate}
\end{ex}

\noindent\textbf{Remarks on log structures on varieties}. It is worth to mention that many other additional structures for generalized pairs have been 
proposed in recent years. For example:
\begin{enumerate}
    \item Hashizume \cite{Has20} proved the non-vanishing and log abundance for projective glc g-pairs $(X,B,\Mm+\bar{A})$ where $A$ is an ample $\Rr$-divisor. This structure is also used in the proof of the cone theorem for glc g-pairs \cite[Theorem 1.3]{HL21}.
    \item The proof of the cone theorem and the existence of flips for glc g-pairs heavily rely on the structure of generalized pairs with ``generically trivial nef part", that is, generalized pairs $(X,B,\Mm)/U$ such that $\Mm_X\sim_{\Rr}0$ over an open subset of $U$ \cite[Theorem 1.1]{HL21}.
\end{enumerate} 
We remark that the structures in both (1) and (2) are special cases of g-pairs with $\bb$-log abundant nef part, but these structures are not preserved under the canonical bundle formula.

We summarize different types of log structures on varieties and their behavior under the non-vanishing conjecture, abundance conjecture, canonical bundle formula, sub-adjunction formula, and the minimal model program in the following table. 

\begin{table}[ht]
   \caption{Different log structures on varieties}
    \label{tab:log structures}
\begin{center}
        \begin{tabular}{|c|c|c|c|c|c|c|}
\hline
& klt & lc & glc(+ample) & glc ($\Mm|_{U^0}\sim_{\Rr}\bf{0}$) & glc, $\Mm$ $\bb$-ab. & glc \\
\hline
Non-vanishing? & C &  C & T & C &  C & F \\
\hline
Abundance? & C &  C & C & C & C & F \\
\hline
Preserved under cbf? & C$^*$ & C & F & F & T & T \\
\hline
Preserved under sub-adj? & C$^*$ & C & T & F & T & T \\
\hline
Preserved under MMP? & T & T & F & T & T & T \\
\hline
\end{tabular}
\end{center}
\begin{flushleft}
C: Conjecturally true. T: True.  F: False. C$^*$: True, but we lose control of coefficients. Only conjecturally true if we want to control the coefficients.
\end{flushleft}
\end{table}

\medskip

\noindent\textbf{Acknowledgement}.  The authors would like to thank Guodu Chen, Christopher D. Hacon, Jingjun Han, Yuchen Liu, Yujie Luo, Fanjun Meng, and Qingyuan Xue for useful discussions. The first and the third author are partially supported by NSF research grants no: DMS-1801851, DMS-1952522 and by a grant from the Simons Foundation; Award Number: 256202.

\section{Preliminaries}

We will freely use the notation and definitions from \cite{KM98,BCHM10}. 

\subsection{Divisors and \texorpdfstring{$\bb$-divisors}{}}
\begin{deflem}[{\cite[III.\S 1]{Nak04}, \cite[Definition-Lemma 3.3.1]{BCHM10}, \cite[Lemma 4.1(1)]{BH14}},{\cite[Section 2.2]{Hu20}}]\label{deflem: Nakayama-Zariski decomposition}
We have the following:
\begin{enumerate}
    \item Let $X$ be a normal projective variety, $C$ a prime divisor over $X$, and $D$ and $\Rr$-Cartier $\Rr$-divisor on $X$ such that $|D|_{\Rr}\not=\emptyset$. We define
$$o_C(D):=\inf\{\mult_CB'\mid B'\in |B|_{\Qq}\}.$$
\item Let $X$ be a smooth projective variety, $D$ a pseudo-effective $\Rr$-divisor on $X$, and $A$ an ample $\Qq$-divisor on $X$. For any prime divisor $C$ on $X$, we define
$$\sigma_C(D):=\lim_{\epsilon\rightarrow 0+}o_C(D+\epsilon A).$$
Then $\sigma_C(D)$ is well-defined and independent of the choice of $A$. Moreover, there are only finitely many prime divisors $C$ such that $\sigma_C(D)>0$. We define $N_{\sigma}(D):=\sum_C\sigma_C(D)C$.
\item Let $X$ be a normal projective variety, $D$ a pseudo-effective $\Rr$-Cartier $\Rr$-divisor on $X$, and $f: Y\rightarrow X$ a resolution of $X$. We define $N_{\sigma}(D):=f_*N_{\sigma}(f^*D)$. Then $N_{\sigma}(D)$ is independent of the choice of $f$.
\end{enumerate}
\end{deflem}

\begin{defn}[$\bb$-divisors]\label{defn: b divisors} Let $X$ be a normal quasi-projective variety. We call $Y$ a \emph{birational model} over $X$ if there exists a projective birational morphism $Y\to X$. 

Let $X\dashrightarrow X'$ be a birational map. For any valuation $\nu$ over $X$, we define $\nu_{X'}$ to be the center of $\nu$ on $X'$. A \emph{$\bb$-divisor} $\Dd$ over $X$ is a formal sum $\Dd=\sum_{\nu} r_{\nu}\nu$ where $\nu$ are valuations over $X$ and $r_{\nu}\in\mathbb R$, such that $\nu_X$ is not a divisor except for finitely many $\nu$. If in addition, $r_{\nu}\in\Qq$ for every $\nu$, then $\Dd$ is called a \emph{$\Qq$-$\bb$-divisor} over $X$. The \emph{trace} of $\Dd$ on $X'$ is the $\Rr$-divisor
$$\Dd_{X'}:=\sum_{\nu_{i,X'}\text{ is a divisor}}r_i\nu_{i,X'}.$$
If $\Dd_{X'}$ is $\Rr$-Cartier and $\Dd_{Y}$ is the pullback of $\Dd_{X'}$ on $Y$ for any birational model $Y$ of $X'$, we say that $\Dd$ \emph{descends} to $X'$ and write $\Dd=\overline{\Dd_{X'}}$.  We let $\bm{0}$ be the $\bb$-divisor $\bar{0}$.

Let $X\rightarrow U$ be a projective morphism and assume that $\Dd$ is a $\bb$-divisor over $X$ such that $\Dd$ descends to some birational model $Y$ over $X$. If $\Dd_Y$ is nef$/U$ (resp. semi-ample$/U$), then we say that $\Dd$ is \emph{nef}$/U$ (resp. semi-ample$/U$). If $\Dd_Y$ is a Cartier divisor, then we say that $\Dd$ is \emph{$\bb$-Cartier}. If $\Dd$ can be written as an $\Rr_{\geq 0}$-linear combination of nef$/U$ $\bb$-Cartier $\bb$-divisors, then we say that $\Dd$ is \emph{NQC}$/U$.
\end{defn}

\subsection{Generalized pairs}

We will follow the original definitions in \cite{BZ16} but will adopt the same notation as in \cite{HL21}. Notice that there are some small differences with the definitions in \cite{HL21}: in this paper, all generalized (sub)-pairs are assumed to be NQC.

\begin{defn}[Generalized pairs]\label{defn: g-pairs}
A \emph{generalized sub-pair} (\emph{g-sub-pair} for short) $(X,B,\Mm)/U$ consists of a normal quasi-projective variety $X$ associated with a projective morphism $X\rightarrow U$, an $\Rr$-divisor $B$ on $X$, and an NQC$/U$ $\bb$-divisor $\Mm$ over $X$, such that $K_X+B+\Mm_X$ is $\Rr$-Cartier. If $B$ is a $\Qq$-divisor and $\Mm$ is a $\Qq$-$\bb$-divisor, then we say that $(X,B,\Mm)/U$ is a \emph{$\Qq$-g-sub-pair}.

If $\Mm=\bm{0}$, a g-sub-pair $(X,B,\Mm)/U$ is called a \emph{sub-pair} and is denoted by $(X,B)$ or $(X,B)/U$. 

If $U=\{pt\}$, we usually drop $U$ and say that $(X,B,\Mm)$ is a \emph{projective}.  If $U$ is not important, we may also drop $U$. This usually happens when we emphasize the structures of $(X,B,\Mm)$ that are independent of the choice of $U$, such as the singularities of $(X,B,\Mm)$. 

A g-sub-pair (resp. $\Qq$-g-sub-pair) $(X,B,\Mm)/U$ is called a \emph{g-pair} (resp. \emph{$\Qq$-g-pair}) if $B\geq 0$. A sub-pair $(X,B)$ is called a \emph{pair} if $B\geq 0$.
\end{defn}

\begin{defn}[Singularities of generalized pairs]\label{defn: sing of g-pairs}
	Let $(X,B,\Mm)/U$ be a g-(sub-)pair. For any prime divisor $E$ and $\mathbb R$-divisor $D$ on $X$, we define $\mult_{E}D$ to be the \emph{multiplicity} of $E$ along $D$.  Let $h:W\to X$
	be any log resolution of $(X,\Supp B)$ such that $\Mm$ descends to $W$, and let
	$$K_W+B_W+\Mm_W:=h^*(K_X+B+\Mm_X).$$
	The \emph{log discrepancy} of a prime divisor $D$ on $W$ with respect to $(X,B,\Mm)$ is $1-\mult_{D}B_W$ and it is denoted by $a(D,X,B,\Mm).$
	
	We say that $(X,B,\Mm)$ is \emph{(sub-)glc} (resp. \emph{(sub-)gklt}) if $a(D,X,B,\Mm)\ge0$ (resp. $>0$) for every log resolution $h: W\to X$ as above and every prime divisor $D$ on $W$. We say that $(X,B,\Mm)$ is \emph{gdlt} if $(X,B,\Mm)$ is glc, and there exists a closed subset $V\subset X$, such that
\begin{enumerate}
    \item $X\backslash V$ is smooth and $B_{X\backslash V}$ is simple normal crossing, and
    \item for any prime divisor $E$ over $X$ such that $a(E,X,B,\Mm)=0$, $\Center_XE\not\subset V$ and $\Center_XE\backslash V$ is an lc center of $(X\backslash V,B|_{X\backslash V})$.
\end{enumerate}
If $\Mm=\bm{0}$ and $(X,B,\Mm)$ is (sub-)glc (resp, (sub-)gklt, gdlt), we say that $(X,B)$ is (sub-)lc (resp. (sub-)klt, dlt).
	    
	 Suppose that $(X,B,\Mm)$ is sub-glc. A \emph{glc place} of $(X,B,\Mm)$ is a prime divisor $E$ over $X$ such that $a(E,X,B,\Mm)=0$. A \emph{glc center} of $(X,B,\Mm)$ is the center of a glc place of $(X,B,\Mm)$ on $X$. The \emph{non-gklt locus} $\Ngklt(X,B,\Mm)$ of $(X,B,\Mm)$ is the union of all glc centers of $(X,B,\Mm)$.
\end{defn}

\begin{deflem}[{cf.~\cite[Proposition 3.9]{HL18}}]\label{deflem: gdlt modification}
Let $(X,B,\Mm)/U$ be a glc g-pair. Then there exists a birational morphism $f: Y\rightarrow X$ and a glc g-pair $(Y,B_Y,\Mm)/U$, such that
\begin{enumerate}
\item $(Y,B_Y,\Mm)$ is $\Qq$-factorial gdlt,
    \item $K_Y+B_Y+\Mm_Y=f^*(K_X+B+\Mm_X)$, and
    \item any $f$-exceptional divisor is a component of $\lfloor B_Y\rfloor$.
\end{enumerate}
For any birational morphism $f$ and $(Y,B_Y,\Mm)$ which satisfies (1-3), $f$ will be called a \emph{gdlt modification} of $(X,B,\Mm)$.
\end{deflem}

\begin{deflem}
Let $(X,B,\Mm)/U$ be a glc g-pair and $S\subset\lfloor B\rfloor$ a prime divisor with normalization $\tilde S$. Then there is a naturally defined g-pair structure $(S,B_S,\Mm^S)/U$ such that $$K_S+B_S+\Mm^S_S=(K_X+B+\Mm_X)|_S$$
(cf. \cite[Definition 4.7]{BZ16}). We say that $(S,B_S,\Mm^S)/U$ is the glc g-pair induced by the adjunction. If $(X,B,\Mm)$ is glc, then $(S,B_S,\Mm^S)$ is glc \cite[Definition 4.7, Remark 4.8]{BZ16}, and if $(X,B,\Mm)$ is gdlt, then $(S,B_S,\Mm^S)$ is gdlt \cite[Lemma 2.6]{HL18}.
\end{deflem}

\begin{defn}\label{defn: b-birational g-pair}
Let $(X,B,\Mm)/U$ and $(X',B',\Mm)/U$ be two g-sub-pairs and $f: X\dashrightarrow X'$ a birational map. We say that $f: (X,B,\Mm)\dashrightarrow (X',B',\Mm)$ is \emph{B-birational} if there exists a common resolution $p: W\rightarrow X$ and $q: W\rightarrow Y$ such that $p^*(K_X+B+\Mm_X)=q^*(K_{X'}+B'+\Mm_{X'})$ and $q=f\circ p$. In addition, if $f$ does not extract any divisor, then we say that  $f: (X,B,\Mm)\dashrightarrow (X',B',\Mm)$  is a \emph{B-birational contraction}.
\end{defn}

\subsection{Canonical bundle formula}

\begin{defn}[Canonical bundle formula]
Let $(X,B,\Mm)/U$ be a g-sub-pair and $f: (X,B,\Mm)\rightarrow Z$ a projective surjective morphism$/U$. A \emph{canonical bundle formula} of $f: (X,B,\Mm)\rightarrow Z$ is a formula of the form
$$K_X+B+\Mm_X\sim_{\Rr}f^*(K_Z+B_Z+\Mm^Z_Z),$$
such that $(Z,B_Z,\Mm^Z)/U$ is a g-sub-pair. We say that $(Z,B_Z,\Mm^Z)/U$ is a \emph{g-sub-pair induced by a canonical bundle formula} of $f: (X,B,\Mm)\rightarrow Z$.
\end{defn}

\begin{defn}[Discrepancy $\bb$-divisor,  cf. {\cite[Example 2.6]{FS20}}]
Let $(X,B,\Mm)/U$ be a g-sub-pair. We define $\bb$-divisors
$\Aa(X,B,\Mm)$ and $\Aa^*(X,B,\Mm)$ in the following way: for any birational morphism $f: Y\rightarrow X$, we define
$$\Aa(X,B,\Mm)_Y:=K_Y+\Mm_Y-f^*(K_X+B+\Mm_X),\text{ and }\Aa^*(X,B,\Mm)_Y:=\Aa(X,B,\Mm)^{>-1}_Y.$$
\end{defn}

\begin{defn}[Glc-trivial fibration, cf. {\cite[Definition 2.19]{FS20}}]
Let $(X,B,\Mm)/U$ be a g-sub-pair and $f: X\rightarrow Z$ a contraction$/U$. If
\begin{enumerate}
    \item $(X,B,\Mm)$ is sub-glc over the generic point of $Z$,
    \item $\rk f_*\mathcal{O}_X(\lceil\Aa^*(X,B,\Mm)\rceil)=1$, and
    \item $K_X+B+\Mm_X\sim_{\Rr,Z}0$,
\end{enumerate}
then we say that $f: (X,B,\Mm)\rightarrow Z$ is a glc-trivial fibration$/U$. Moreover, if $\Mm=\bm{0}$, then we say that $f: (X,B)\rightarrow Z$ is an lc-trivial fibration$/U$. 

If we additionally have
\begin{enumerate}
    \item[(1')] $(X,B,\Mm)$ is sub-gklt over the generic point of $Z$,
\end{enumerate}
then we say that $f: (X,B,\Mm)\rightarrow Z$ is a gklt-trivial fibration$/U$. Moreover, if $\Mm=\bm{0}$, then we say that $f: (X,B)\rightarrow Z$ is a klt-trivial fibration$/U$. 
\end{defn}

\begin{defn}[Glc-trivial morphism]
Let $(X,B,\Mm)/U$ be a sub-glc g-sub-pair and $f: X\rightarrow Z$ a projective surjective morphism between normal quasi-projective varieties. Let $X\xrightarrow{\tau}\tilde Z\xrightarrow{\gamma}Z$ be the Stein factorization of $f$. We say that $f: (X,B,\Mm)\rightarrow Z$ is a \emph{glc-trivial morphism}$/U$ if $K_X+B+\Mm_X\sim_{\Rr,Z}0$ and $\tau: (X,B,\Mm)\rightarrow\tilde Z$ is a glc-trivial fibration$/U$.
\end{defn}

\begin{defn}
Let $f: X\rightarrow Z$ a projective surjective morphism between normal quasi-projective varieties, $(X,B,\Mm)/U$ a g-sub-pair that is sub-glc over the generic point of $Z$, such that $K_X+B+\Mm_X\sim_{\Rr,Z}0$. Let $D$ be a prime divisor over $Z$, and let $h_Z: Z'\rightarrow Z$ be a birational morphism such that $D$ is on $Z'$. We let $X'$ be the normalization of the main component of $X\times_ZZ'$, $f': X'\rightarrow Z'$ and $h: X'\rightarrow X$ the induced morphisms, and $K_{X'}+B'+\Mm_{X'}:=h^*(K_X+B+\Mm_X)$. We define
$$t_D(X,B,\Mm;f):=\sup\{(X',B'+tf'^*D,\Mm)\text{ is sub-glc over the generic point of } D\}.$$
Note that $f^*D$ is well-defined over the generic point of $D$ even if $D$ is not $\Qq$-Cartier. It is clear that $t_D(X,B,\Mm;f)$ is independent of the choice of $h_Z$.
\end{defn}

\subsection{Canonical bundle formulas}\label{subsection: scbf}

\begin{deflem}\label{deflem: scbf}
Let $f: (X,B,\Mm)\rightarrow Z$ be a glc-trivial fibration. A \emph{(Kodaira type) canonical bundle formula} for $f: (X,B,\Mm)\rightarrow Z$ is a canonical bundle formula
$$K_X+B+\Mm_X\sim_{\Rr}f^*(K_Z+B_Z+\Mm^Z_Z)$$
such that for any prime divisor $D$ over $Z$ and any birational morphism $h_Z: Z'\rightarrow Z$ such that $D$ is on $Z'$, we have $\mult_{D}B_{Z'}=1-t_D(X,B,\Mm;f)$ where $K_{Z'}+B_{Z'}+\Mm^Z_{Z'}:=h_Z^*(K_Z+B_Z+\Mm^Z_Z)$.

If  $(Z,B_Z,\Mm^Z)/U$ is a g-sub-pair, then we call it a \emph{g-sub-pair induced by (a canonical bundle formula of)} $f: (X,B,\Mm)\rightarrow Z$. Moreover, by our construction, if $(X,B,\Mm)$ is sub-glc (resp. sub-gklt, glc, gklt), then $(Z,B_Z,\Mm^Z)$ is sub-glc (resp. sub-gklt, glc, gklt). Note that we cannot find any reference answering whether $(Z,B_Z,\Mm^Z)/U$ is a g-sub-pair in general. Nevertheless, by \cite[Chapter 6, Theorem 7]{Fil19} and \cite[Theorem 2.20]{FS20}, if $f: (X,B,\Mm)\rightarrow Z$ a glc-trivial 
fibration$/U$ such that $(X,B,\Mm)$ is a $\Qq$-g-sub-pair and
\begin{enumerate}
    \item either $B\geq 0$ over the generic point of $Z$, or
    \item $\Mm$ is semi-ample$/Z$,
\end{enumerate}
then $\Mm^Z$ is always nef$/U$, and in these cases, $(Z,B_Z,\Mm^Z)/U$ will be guaranteed as a g-sub-pair. We will show later that the $\Qq$-coefficient assumption can be removed (Theorem \ref{thm: cbf for R-coefficient}).
\end{deflem}

\begin{deflem}\label{deflem: cbf finite}
Let $(X,B,\Mm)/U$ be a sub-glc g-sub-pair and $f: (X,B,\Mm)\rightarrow Z$ a finite morphism$/U$ such that $K_X+B+\Mm_X\sim_{\Rr,Z}0$. Let $Z^o$ be the smooth locus of $Z$ and $X^o:=f^{-1}(Z^o)$. By Hurwitz formula, we have
$$K_{X^o}=(f|_{X^o})^*K_{Z^o}+R^o,$$
where $R^o$ is the ramification divisor of $f|_{X^o}$. We let $R$ be the closure of $R^o$ in $X^o$. We define $$B_Z:=\frac{1}{\deg f}f_*(R+B).$$ Now consider any proper birational morphism $h_Z: Z'\rightarrow Z$. We let $X'$ be the normalization of the main component of $Z'\times_{Z}X$, and $f':X'\rightarrow Z'$ the induced morphism. We define
$$\Mm^Z_{Z'}:=\frac{1}{\deg f}f'_*\Mm_{X},$$
and let $\Mm^Z$ be the corresponding $\bb$-divisor. By \cite[Theorem 4.5]{HL19},  $(Z,B_Z,\Mm^Z)/U$ is a sub-glc g-sub-pair, and we will call
$$K_X+B+\Mm_X\sim_{\Rr}f^*(K_Z+B_Z+\Mm^Z_Z).$$
the \emph{(Fujino-Gongyo type) canonical bundle formula} of $f: (X,B,\Mm)\rightarrow Z$. We will call $(Z,B_Z,\Mm)/U$  a \emph{g-sub-pair induced by  (a canonical bundle formula of)} $f: (X,B,\Mm)\rightarrow Z$.
\end{deflem}

\begin{deflem}\label{deflem: cons g-pair finite morphism}
Let $(X,B,\Mm)/U$ be a sub-glc g-sub-pair and $f: (X,B,\Mm)\rightarrow Z$ a glc-trivial morphism$/U$.  Let $X\xrightarrow{\tau}\tilde Z\xrightarrow{\gamma}Z$ be the Stein factorization of $f$. Assume that we have a g-sub-pair $(\tilde Z,B_{\tilde Z},\Mm^{\tilde Z})/U$ induced by  $\tau: (X,B,\Mm)\rightarrow \tilde Z$, and we let $(Z,B_Z,\Mm^Z)/U$ be the g-sub-pair induced by $\gamma: (\tilde Z,B_{\tilde Z},\Mm^{\tilde Z})\rightarrow Z$. We say that $(Z,B_Z,\Mm^Z)/U$ is a \emph{g-sub-pair induced by (a canonical bundle formula of)} $f: (X,B,\Mm)\rightarrow Z$.

Moreover, $f: (X,B,\Mm)\rightarrow Z$ is called \emph{good} if $\gamma^*\Mm^Z=\Mm^{\tilde Z}$, and $f: (X,B,\Mm)\rightarrow Z$ is called a \emph{gdlt-trivial morphism$/U$} if $(Z,B_Z,\Mm^Z)$ and $(X,B,\Mm)$ are both gdlt. We remark that if $\Mm=\bf{0}$, then a gdlt-trivial morphism$/U$ is called a good dlt model in \cite{Hu20}.
\end{deflem}

\subsection{Numerical properties}

\begin{defn}
Let $\pi: X\rightarrow U$ be a projective morphism from a normal variety to a variety and $D$ an $\Rr$-divisor on $X$. We define $\kappa(X/U,D),\kappa_{\iota}(X/U,D),\kappa_{\sigma}(X/U,D)$ to be the \emph{Iitaka dimension}, \emph{invariant Iitaka dimension}, and \emph{numerical Iitaka dimension} of $D$ over $U$ respectively. For formal definitions and their basic properties, we refer the reader to \cite[Chapters II,V]{Nak04}, \cite[Section 2]{Cho08}, and \cite[Section 2]{HH20}.
\end{defn}

\begin{defn}\label{defn: b abundance}
Let $X\rightarrow U$ be a projective morphism between normal varieties, $D$ an $\Rr$-divisor on $X$, $\Dd$ a $\bb$-divisor over $X$, and $Y$ a birational model of $X$ such that $\Dd$ descends to $Y$. We say that $D$ is \emph{abundant}$/U$ if $\kappa_{\iota}(X/U,D)=\kappa_{\sigma}(X/U,D)$, and we say that $\Dd$ is $\bb$-\emph{abundant}$/U$ if $\Dd_Y$ is abundant$/U$.
\end{defn}

\begin{defn}\label{defn: b log abundance}
Let $(X,B,\Mm)/U$ be a sub-glc g-sub-pair, $D$ an $\Rr$-divisor on $X$, $\Dd$ a $\bb$-divisor over $X$, $Y$ a birational model of $X$ with induced birational morphism $f: Y\rightarrow X$ such that $\Dd$ descends to $Y$, and $(Y,B_Y,\Mm)/U$ the sub-glc g-sub-pair defined by $K_Y+B_Y+\Mm_Y:=f^*(K_X+B+\Mm_X)$. We say that $D$ is \emph{log abundant}$/U$ with respect to $(X,B,\Mm)$ if $D$ is log abundant$/U$, and for any glc center $W$ of $(X,B,\Mm)$ with normalization $W^{\nu}$, $D|_{W^{\nu}}$ is abundant$/U$. We say that $\Dd$ is $\bb$\emph{-log abundant}$/U$ with respect to $(X,B,\Mm)$ if $D_Y$ is log abundant$/U$ with respect to $(Y,B_Y,\Mm)$.
\end{defn}

\begin{lem}\label{lem: limit of nef and abundant divisors is nef and abundant}
Let $X\rightarrow U$ be a projective morphism from a normal variety to a variety, $\{D_i\}_{i=1}^{+\infty}$ a sequence of $\Rr$-divisors on $X$, and $D$ an $\Rr$-divisor on $X$, such that
\begin{enumerate}
    \item $D\geq D_i$ for each $i$,
    \item $\lim_{i\rightarrow+\infty}||D-D_i||=0$, and
    \item each $D_i$ is abundant$/U$.
\end{enumerate}
Then $D$ is abundant$/U$. Moreover, if each $D_i$ is nef$/U$, then $D$ is nef$/U$.
\end{lem}
\begin{proof}
Possibly passing to a subsequence, we may assume that $D_2\geq D_1$ and $\Supp(D_2-D_1)=\Supp(D-D_1)$. Thus $D_2\geq\lambda D_1+(1-\lambda)D$ for some $\lambda\in (0,1)$. Since $D_2$ is nef and abundant, $\kappa_{\iota}(X/U,D_2)=\kappa_{\sigma}(X/U,D_2)\geq 0$. Since $D\geq D_2$, $\kappa_{\iota}(X/U,D_2)\geq 0$, and $\kappa_{\sigma}(X/U,D_2)\geq 0$. Then
$$\kappa_{\iota}(X/U,D)\geq\kappa_{\iota}(X/U,D_2)=\kappa_{\sigma}(X/U,D_2)\geq\kappa_{\sigma}(X/U,D).$$
So $D$ is abundant$/U$.  The moreover part is obvious.
\end{proof}

\begin{lem}[{\cite[Lemma 2.8]{Hu20}}]\label{lem: equiv def for nef and abundant}
Let $\pi: X\to U$ be a projective morphism from a normal variety to a variety, and $D$ a nef$/U$ $\Rr$-divisor on $X$. Then the following conditions are equivalent:
\begin{enumerate}
    \item $D$ is abundant$/U$.
    \item There exists a birational model $h: X'\to X$ and a surjective morphism$/U$ $g: X'\to Y$ such that $g^*H\sim_{\Rr}h^*D$, where $H$ is a nef and big$/U$ $\Rr$-divisor on $Y$.
    \item Given a sufficiently general fiber $F$ of $\pi$, $o_\Gamma(D|_F)=0$ for every prime divisor $\Gamma$ over $F$.
\end{enumerate}
\end{lem}

\subsection{On real coefficients}

%\begin{lem}\label{lem: easy perturbation lemma}
%Let  $r_1,\dots,r_c$ be real numbers such that $1,r_1,\dots,r_c$ are linearly independent over $\Qq$, and $s_1,s_2:\mathbb R^{c+1}\rightarrow\mathbb R$ $\Qq$-linear functions. Assume that $s_1(1,\bm{r})=qs_2(1,\bm{r})$ for some $q\in\mathbb Q$. Then $s_1=qs_2$.
%\end{lem}
%\begin{proof}
%We have $s_1(x_0,\dots,x_c)=\sum_{i=0}^cm_ix_i$ and $s_2(x_0,\dots,x_c)=\sum_{i=0}^cn_ix_i$ for any $x_0,\dots,x_c$, where $m_0,\dots,m_c,n_0,\dots,n_c\in\mathbb Q$. Then $(m_i-qn_i)+\sum_{i=1}^c(m_i-qn_i)r_i=0$. Since  $1,r_1,\dots,r_c$ are linearly independent over $\Qq$, $m_i=qn_i$ for each $i$. Thus $s_1=qs_2$.
%\end{proof}

\begin{lem}\label{lem: abundant is open in minimal Q affine V}
Let $M$ be an abundant $\Rr$-divisor on $X$ and $V$ the minimal rational affine subspace of $\Weil_\Rr(X)$ which contains $M$. Then there is a rational polytope $P\subset V$ containing $M$ such that any $M'\in P$ is abundant.
\end{lem}
\begin{proof}
Write $M=\sum_{i=1}^n a_iM_i+\sum_{i=1}^nb_i\text{div}(\phi_i)=\sum_{i=1}^{l}c_iE_i$, where $M_i,E_i$ are $\Qq$-divisors, $E_i\ge 0$ and $c_i>0$. Then we have a linear map 
$$
\Phi: \Rr^m\times\Rr^n\times\Rr^l\to \Weil_\Qq(X)\otimes_\Qq\Rr
$$
sending $(\textbf{x},\bf{y},\bf{z})$ to $\sum_{i=1}^n x_iM_i+\sum_{i=1}^ny_i\text{div}(\phi_i)-\sum_{i=1}^{l}z_iE_i$. We note that $\Phi$ is defined over $\Qq$. Then the statement follows easily by considering the $\Qq$-affine space $\Phi^{-1}(0)$ and the following fact:
$$
\sum_{i=1}^ld_iE_i \text{ is abundant iff } \sum_{i=1}^le_iE_i \text{ is abundant where } e_i,d_i>0. 
$$

\end{proof}

\begin{lem}\label{lem: lct is linear}
Let $(X,B)$ be a sub-lc sub-pair and $V$ the minimal rational affine subspace of $\Weil_\Rr(X)$ such that $B\in V$. Let $D\ge0$ be a $\Qq$-Cartier $\Qq$-divisor on $X$ and $$t:=\sup\{s\in\Rr| (X,B+sD)\text{ is sub-lc}\}.$$ Then there exists a rational polytope $P\subset V$ containing $B$, such that for any $B'\in P$,
\begin{enumerate}
    \item $(X,B')$ is sub-lc,
    \item if $a(E,X,B+tD)=0$, then $a(E,X,B'+t'D)=0$ , where $$t':=\sup\{s\in\Rr| (X,B'+sD)\text{ is sub-lc}\}.$$
\end{enumerate}

\end{lem}

\begin{proof}
The minimality of $V$ implies that $K_X+B'$ is $\Qq$-Cartier for any $B'\in V$.
Possibly passing to a log resolution of $(X,\Supp B\cup\Supp D)$, we may assume that $(X,\Supp B\cup\Supp D)$ is log smooth.
%We can reduce to the case that $(X,B+D)$ is snc by considering the canonical maps $f_*:\Weil(Y)\to\Weil(X)$ and $f^*:\text{CDiv}(X)\to\text{CDiv}(Y)$, where $f:Y\to X$ is some log resolution. Thanks to the minimality of $V$ and the snc condition, 
Let $E_1,\dots,E_n$ be the components of $D$ for some positive integer $n$. For each $i$, we have an affine function $t_i: V\to \Qq$ given by $t_i(B'):=\frac{1-\mult_{E_i}B'}{\mult_{E_i}D}$.
Notice that if $(X,B')$ is sub-lc, $a(E_i,X,B'+t'D)=0$ is equivalent to say that $t_i(B')\ge t_j(B')$ for all $1\le j\le n$. Thus $\{B'\in V|~B'\text{ satisfies (1) and (2)}\}$ is a rational polyhedron and we are done.
\end{proof}

\begin{thm}[Kodaira-type formula for $\Rr$-coefficients]\label{thm: cbf for R-coefficient}
Let $(X,B,\Mm)/U$ be a g-sub-pair and $f: (X,B,\Mm)\rightarrow Z$ a glc-trivial fibration$/U$, such that
\begin{itemize}
    \item either $B\geq 0$ over the generic point of $Z$, or
    \item $\Mm$ is semi-ample$/Z$ (in particular, this includes the case when $\Mm=\bf{0}$).
\end{itemize}
Then there exists a g-sub-pair $(Z,B_Z,\Mm^Z)$ on $Z$ such that $a(D,Z,B_Z,\Mm^Z)=1-t_D(X,B,\Mm;f)$ for any prime divisor $D$ over $Z$ and $K_X+B+\Mm_X\sim_{\Rr}f^*(K_Z+B_Z+\Mm^Z_Z)$.
\end{thm}
\begin{proof}
We only need to show that $\Mm^Z$ is NQC$/U$. Possibly by replacing $B$ with $B-f^*H$ for some ample divisor on $Z$ we can assume $(X,B,\Mm)$ is sub-glc. By \cite[Lemma 4.1]{HL19}, we may find a positive integer $k$, real numbers $a_1,\dots,a_k$ and $\Qq$-g-sub-pairs $(X,B_i,\Mm_i)/U$,such that $\sum_{i=1}^ka_i=1$, $\sum_{i=1}^ka_iB_i=B$, $\sum_{i=1}^ka_i\Mm_i=\Mm$, $K_X+B_i+\Mm_{i,X}\sim_{\Qq,Z}0$, and either each $B_i\geq 0$ over the generic point of $Z$ or $\Mm_i$ is semi-ample$/Z$. By \cite[Chapter 6, Theorem 7]{Fil19} and \cite[Theorem 2.20]{FS20}, for every $1\leq i\leq c+1$, there exists a $\Qq$-g-sub-pair $(Z,B_{Z,i},\Mm^{Z,i})/U$ such that
$$K_X+B_i+\Mm_{i,X}\sim_{\Qq}f^*(K_{Z}+B_{Z,i}+\Mm^{Z,i}_{Z}),$$
and $a(D,Z,B_{Z,i},\Mm^Z_i)=1-t_D(X,B_i,\Mm_i;f)$ for any prime divisor $D$ over $Z$. We only need to choose suitable $(X,B_i,\Mm_i)$ such that  $\sum_{i=1}^ka_it_D(X,B_i,\Mm_i;f)=t_D(X,B,\Mm;f)$ for any prime divisor $D$ over $Z$. By applying the weak-semistable reduction \cite{AK00} (see \cite[Theorem B.6]{Hu20} for details), we only need to consider finitely many prime divisors $D$ over $Z$. The theorem follows from Lemma \ref{lem: lct is linear}, \cite[Lemma 3.3]{Hu20}, and \cite[Lemma 4.1]{HL19} after replacing $(X,B_i,\Mm_i)$ and $a_i$ for each $i$ so that $||B_i-B||$ is sufficiently small and $B_i$ is contained in the minimal affine subspace of $\Weil_{\Qq}(X)$ which contains $B$.
\end{proof}

%\begin{lem}
%Let $(X,B,\Mm)/U$ be a g-sub-pair and $f:(X,B,\Mm)/U\to Z$ is a glc-trivial fibration$/U$. Then there are finitely many positive real numbers $c_\lambda$ satisfying $\sum_\lambda c_\lambda=1$ and  $\QQ$-g-sub-pairs $(X,B^\lambda,\Mm^\lambda)$ such that $K_X+B+\Mm=\sum_\lambda c_\lambda(K_X+B^\lambda+\Mm^\lambda_X)$ and $f:(X,B^\lambda,\Mm^\lambda)$ is a glc-trivial fibration for each $\lambda$
%\end{lem}

\section{Canonical bundle formula: the klt case}

In this section, we prove Theorem \ref{thm: log abundance gcbf} when $(X,B,\Mm)$ is klt. More precisely, we have the following theorem:

\begin{thm}\label{thm: abundant gcbf klt case}
Let $(X,B,\Mm)/U$ be a $\Qq$-g-sub-pair such that $\Mm$ is $\bb$-abundant$/U$, and $f: (X,B,\Mm)\rightarrow Z$ a gklt-trivial fibration$/U$. Let
$$K_X+B+\Mm_X\sim_{\Qq}f^*(K_Z+B_Z+\Mm^Z_Z)$$
be a canonical bundle formula. Then $\Mm^Z$ is $\bb$-abundant$/U$.
\end{thm}

\begin{lem}\label{lem: abundant gcbf klt ssample case}
Let $(X,B,\Mm)/U$ be a $\Qq$-g-sub-pair such that $\Mm=\overline{\Mm_X}$, and $f: (X,B,\Mm)\rightarrow Z$ a gklt-trivial fibration$/U$, such that $\Mm$ is semi-ample$/U$. Let
$$K_X+B+\Mm_X\sim_{\Qq}f^*(K_Z+B_Z+\Mm^Z_Z)$$
be a canonical bundle formula. Then $\Mm^Z$ is $\bb$-abundant$/U$.
\end{lem}
\begin{proof}
Possibly replacing $Z$ and $X$ with higher models, we may assume that $\Mm^Z$ descends to $Z$. By \cite[Chapter 6, Theorem 7]{Fil19}, $\Mm^Z$ is nef$/U$, hence $\Mm^Z_Z$ is nef$/U$.

Fix a prime divisor $\Gamma$ over a sufficiently general fiber $F$ of $Z\rightarrow U$. In the following, we will show that $o_{\Gamma}(\Mm^Z_Z|_F)=0$, and the theorem will follow from Lemma \ref{lem: equiv def for nef and abundant}.

Take a resolution $h_Z: Z'\rightarrow Z$ and a prime divisor $D$ on $Z'$, such that $\Gamma$ is on a sufficiently general fiber $F'$ of the induced morphism $Z'\rightarrow U$. Possibly replacing $Z$ with $Z'$ and replacing $X$ accordingly, we may assume that $\Gamma$ is on $F$. Let $F_X$ be the fiber of $X\rightarrow U$ which dominates $F$. By the definition of relative Iitaka dimensions and numerical dimensions, possibly replacing $X$ with $F_X$, $Z$ with $F$, and $U$ with $\{pt\}$, we may assume that $X$ and $Z$ are projective, $U$ is a point, and $\Gamma$ is on $Z$.

Since $\Mm_X$ is semi-ample, we may pick a positive integer $n$ such that $n\Mm_X$ is base-point-free. In particular, we may pick $H\in |2n\Mm_X|$ such that $f_*H$ does not contain the generic point of $\Gamma$ and $(X,B+\frac{1}{2}H)$ is sub-klt over the generic point of $Z$.

\begin{claim}\label{claim: add semiample divisor still klt trivial}
 $f: (X,B+\frac{1}{2}H)\rightarrow Z$ is a klt-trivial fibration$/U$.
\end{claim}
\begin{proof}
Since $(X,B+H)$ is sub-klt, $\rk f_*\mathcal{O}_X(\lceil\mathbf{A}^*(X,B+H)\rceil)>0$. Since
$$a(D,X,B+H)\leq a(D,X,B)=a(D,X,B,\Mm)$$
for any prime divisor $D$ over $X$, 
$$\rk f_*\mathcal{O}_X(\lceil\mathbf{A}^*(X,B+H)\rceil)\leq\rk f_*\mathcal{O}_X(\lceil\mathbf{A}^*(X,B,\Mm)\rceil)=1,$$
hence $\rk f_*\mathcal{O}_X(\lceil\mathbf{A}^*(X,B+H)\rceil)=1$ and $f$ is a klt-trival fibration$/U$.
\end{proof}
\noindent\textit{Proof of Lemma \ref{lem: abundant gcbf klt ssample case} continued}. By Claim \ref{claim: add semiample divisor still klt trivial} and \cite[Theorem 3.3]{Amb05}, we have a canonical bundle formula
$$K_X+B+\frac{1}{2}H\sim_{\Qq}f^*(K_Z+B^H_Z+\Mm^{Z,H}_Z)$$
of $f: (X,B+\frac{1}{2}H)\rightarrow Z$ which induces a gklt g-pair $(Z,B^H_Z,\Mm^{Z,H})/U$, such that $\Mm^{Z,H}$ is $\bb$-abundant. Since $f_*H$ does not contain the generic point of $\Gamma$, $$\lct(X,B+\frac{1}{2}{H};f^*\Gamma)=\lct(X,B;f^*\Gamma)=\glct(X,B,\Mm;f^*\Gamma)$$
over the generic point of $\Gamma$, hence $\mult_{\Gamma}B_Z=\mult_{\Gamma}B^H_Z$. 

Let $h_Z: Z'\rightarrow Z$ be a birational morphism such that $\Mm^{Z,H}$ descends to $Z'$. Let $K_{Z'}+B_{Z'}^H+\Mm^{Z,H}_{Z'}:=h_Z^*(K_Z+B^H_Z+\Mm^{Z,H}_Z)$, $K_{Z'}+B_{Z'}+\Mm^{Z}_{Z'}:=h_Z^*(K_Z+B_Z+\Mm^{Z}_Z)$, and $\Gamma':=(h_Z^{-1})_*\Gamma$. By our construction, $B^H_{Z'}\geq B_{Z'}$. Since $\Mm^{Z,H}_{Z'}$ is abundant,
$$o_{\Gamma'}(\Mm^{Z}_{Z'})=o_{\Gamma'}(\Mm^{Z,H}_{Z'}+B^H_{Z'}-B_{Z'})\leq o_{\Gamma'}(\Mm^{Z,H}_{Z'})=0,$$
hence $o_{\Gamma}(\Mm^Z_Z)=0$. This concludes the proof.
\end{proof}

\begin{proof}[Proof of Theorem \ref{thm: abundant gcbf klt case}]
Since $\Mm$ is $\bb$-abundant, by Lemma \ref{lem: equiv def for nef and abundant}, possibly replacing $(X,B,\Mm)$ with a high resolution, we may assume that $(X,\Supp B)$ is log smooth, $\Mm=\overline{\Mm_X}$, and there exist a surjective morphism$/U$ $g: X\rightarrow Y$ and a big and nef$/U$ $\Qq$-divisor $H$ on $Y$, such that $g^*H=\Mm_X$. Thus there exists a $\Qq$-divisor $E\geq 0$, such that for any $n\gg 0$, there exists an ample$/U$ $\Qq$-divisor $A_n$ on $Y$ such that $H=A_n+\frac{1}{n}E$. We let $F:=g^*E$ and let $\NN_n:=\overline{g^*A_n}$. 

Since $(X,B,\Mm)$ is sub-gklt, $(X,B+\frac{1}{n}F,\NN_n)$ is sub-gklt for any $n\gg 0$. Moreover, since $\Mm=\overline{\Mm_X}$, for any prime divisor $D$ over $X$ and any positive integer $n$,
$$a(E,X,B,\Mm)=a(E,X,B)\geq a(E,X,B+\frac{1}{n}F)=a(E,X,B+\frac{1}{n}F,\NN_n).$$
Thus
$$1=\rk f_*\mathcal{O}_X(\lceil\Aa^*(X,B,\Mm)\rceil)\geq\rk f_*\mathcal{O}_X(\lceil\Aa^*(X,B+\frac{1}{n}F,\NN_n)\rceil)>0$$
for any $n\gg 0$. Therefore, $\rk f_*\mathcal{O}_X(\lceil\Aa^*(X,B+\frac{1}{n}F,\NN_n)\rceil)=1$, hence $f: (X,B+\frac{1}{n}F,\NN_n)\rightarrow Z$ is a gklt-trivial fibration$/U$ for any $n\gg 0$. We let 
$$K_X+(B+\frac{1}{n}F)+\NN_{n,X}\sim_{\Qq}f^*(K_Z+B_{Z,n}+\Mm^{Z,n}_{Z})$$
be a canonical bundle formula for $f: (X,B+\frac{1}{n}F,\NN_n)\rightarrow Z$. Since $\NN_{n}$ descends to $X$ for each $n$ and since $\Supp(B+\frac{1}{n}F)$ is a fixed divisor for $n\gg 0$, possibly replacing $X,Z$ with higher resolutions, we may assume that $\Mm^{Z,n}$ descends to $Z$ for any $n\gg 0$ and $\Mm^Z$ descends to $Z$. 

By the construction of the canonical bundle formulas, we have $B_{Z,n}\geq B_Z$ for any $n$ and $\lim_{n\rightarrow+\infty}||B_{Z,n}-B_Z||=0$. Since $B_{Z,n}+\Mm^{Z,n}_Z\sim_{\Qq}B_Z+\Mm^Z_Z$ for any $n$, possibly replacing $\Mm^{Z,n}$ with $\overline{\Mm^Z_Z+B_Z-B_{Z,n}}$, we may assume that $\Mm^Z_Z\geq\Mm^{Z,n}_Z$ and $\lim_{n\rightarrow+\infty}||\Mm^{Z}_Z-\Mm^{Z,n}_Z||=0$. By Lemma \ref{lem: abundant gcbf klt ssample case}, $\Mm^{Z,n}$ is $\bb$-abundant$/U$ for any $n\gg 0$, hence $\Mm^{Z,n}_Z$ is abundant$/U$ for any $n\gg 0$. By Lemma \ref{lem: limit of nef and abundant divisors is nef and abundant}, $\Mm^Z_Z$ is abundant$/U$, hence $\Mm^Z$ is $\bb$-abundant$/U$.
\end{proof}

\begin{comment}
\begin{proof}[Proof of Corollary \ref{cor: gklt abundance gcbf}]
The $\Qq$-coefficient case follows from Theorem \ref{thm: abundant gcbf klt case}, \cite[Theorem 2.20]{FS20}, and \cite[Theorem 1.3]{Hu20}. The $\Rr$-coefficient case follows from the $\Qq$-coefficient case and the standard theory of Shokurov-type rational polytopes (cf. \cite[Proposition 3.16]{HL18}) and uniform rational polytopes (cf. \cite[Lemma 5.3]{HLS19}, \cite[Therem 1.4]{Che20}).
\end{proof}
\end{comment}

\section{Canonical bundle formula: the general case}

In this section, we prove Theorems \ref{thm: log abundance gcbf} and \ref{thm: sub-adjunction log abundant}. The proof of Theorem \ref{thm: log abundance gcbf} is similar to \cite[Subsection 3.4]{Hu20} (see also \cite[Proposition 4.4]{FL19}) whose ideas were originated in \cite{FG14}. For the reader's convenience, we give a full proof here.

First, we set up the following condition, which will be applied several times in this section.
\begin{cond}\label{cond: weak-semistabilized}
$(X,B,\Mm)/U,(Z,B_Z,\Mm^Z)/U$, $f: X\rightarrow Z$, and $\Delta$ and $\Delta_Z$ are given as follows:
\begin{enumerate}
\item $(X,B,\Mm)/U$ and $(Z,B_Z,\Mm^Z)/U$ are sub-glc $\Qq$-g-sub-pairs.
\item $\Delta,\Delta_Z$ are reduced divisors, $(Z,\Delta_Z)$ is log smooth, $\Supp B\subset\Delta$, and $\Supp B_Z\subset\Delta_Z$.
\item $\Mm$ descends to $X$ and $\Mm^Z$ descends to $Z$.
\item $f: (X,\Delta)\rightarrow (Z,\Delta_Z)$ is a toroidal from a quasi-smooth toroidal variety to a smooth variety, and $f^{-1}(Z\backslash\Delta_Z)$ is smooth.
\item $f: (X,B,\Mm)\rightarrow Z$ is a glc-trivial fibration$/U$.
\item $$K_X+B+\Mm_X\sim_{\Qq}f^*(K_Z+B_Z+\Mm^Z_Z)$$
is a canonical bundle formula.
\item All fibers of $f$ are reduced, and $f$ is flat. In particular, $f$ is equi-dimensional.
\end{enumerate}
\end{cond}

The statement and proof of the following lemma is similar to \cite[Lemma 3.33]{Hu20}.
\begin{lem}\label{lem: hu20 3.33 gpair}
Suppose that $(X,B,\Mm)/U,(Z,B_Z,\Mm^Z)/U$, $f: X\rightarrow Z$, and $\Delta$ and $\Delta_Z$ satisfy Condition \ref{cond: weak-semistabilized}. Let $D$ be a prime divisor on $Z$ such that $D\not\subset\Delta_Z$.
\begin{enumerate}
    \item  If $(Z,D+\Delta_Z)$ is log smooth, then $t_D(X,B,\Mm;f)=1$ (see Definition-Lemma \ref{deflem: scbf}) and $(X,B+f^*D,\Mm)$ is sub-glc.
    \item Let $S\subset B^{=1}$ be a prime divisor that is horizontal$/Z$, and let $(S,B_S,\Mm^S)/U$ be the g-pair induced by the adjunction $K_S+B_S:=(K_X+B)|_S$ and $\Mm^S:=\overline{\Mm_X|_S}$, and $f|_S: S\rightarrow Z$ the induced morphism. Then $t_D(S,B_S,\Mm^S;f|_S)=1$.
    \item Let $S\subset B^{=1}$ be a prime divisor that is vertical$/Z$ such that $T:=f(S)\subset Z$ is a prime divisor. Let $(S,B_S,\Mm^S)/U$ be the g-pair induced by the adjunction $K_S+B_S:=(K_X+B)|_S$ and $\Mm^S:=\overline{\Mm_X|_S}$, and $f|_S: S\rightarrow T$ the induced morphism. Let $D_T$ be a prime divisor on $T$ such that $D_T\not\subset(\Delta_Z-T)|_T$. Then $t_{D_T}(S,B_S,\Mm^S;f|_S)=1$.
\end{enumerate}
\end{lem}
\begin{proof}
\begin{enumerate}
    \item Let $t:=t_D(X,B,\Mm;f)$, then by our assumption, $t>0$. Thus there exists a glc center $S$ of $(X,B+tf^*D,\Mm)$ supported on $f^*D$ which is not a glc center of $(X,B,\Mm)$. Therefore, there exists a glc center $T$ of $(Z,B_Z+tD,\Mm^Z)$ such that $T\subset D$. Notice that $(Z,D+\Delta_Z)$ is log smooth and $\Mm^Z$ descends on $Z$, we must have $t=1$ since $D\not\subset B_Z$ and thus $(X,B+f^*D,\Mm)$ is sub-glc.
    \item Possibly shrinking $Z$ to a neighborhood of the generic point of $D$ and shrinking $X$ accordingly, we may assume that $\Delta_Z=0$. (2) follows from (1) and a direct adjunction computation.
    \item Possibly shrinking $Z$ and $X$, we may assume that $D_T$ is smooth and $\Delta_Z=T$. By \cite[Lemma 3.32]{Hu20}, there exists a prime divisor $D\subset Z$ such that $D_T\subset D$ and $(Z,D+\Delta_Z)$ is log smooth. By (1), $(X,B+f^*D,\Mm)$ is sub-glc, hence $t_D(X,B,\Mm;f)=t_{D_T}(S,B_S,\Mm^S;f|_S)=1.$
\end{enumerate}
\end{proof}

The statement and the proof of the following lemma is similar to \cite[Lemma 3.35]{Hu20}. It is also similar to \cite[Proposition 4.4]{FL19} and \cite[Proof of Theorem 1.1]{FG14}.

\begin{lem}\label{lem: hu20 3.35 gpair}
Suppose that $(X,B,\Mm)/U,(Z,B_Z,\Mm^Z)/U$, $f: X\rightarrow Z$, $\Delta$ and $\Delta_Z$ satisfy Condition \ref{cond: weak-semistabilized}. Let $S\subset B^{=1}$ be a prime divisor and $F$ a general fiber of $f$, such that
\begin{enumerate}
    \item $T:=f(S)$ is either $Z$ or a prime divisor on $Z$ and
    \item $-B^{<0}|_F=N_{\sigma}((K_X+B^{>0}+\Mm_X)|_F)$.
\end{enumerate}
Let $(S,B_S,\Mm^S)/U$ be the sub-glc g-sub-pair induced by the adjunction $K_S+B_S:=(K_X+B)|_S$ and $\Mm^S:=\overline{\Mm_X|_S}$, and $f|_S: S\rightarrow T$ the induced morphism. Then there exists a canonical bundle formula
$$K_S+B_S+\Mm^S_S\sim_{\Qq}(f|_S)^*(K_T+B_T+\Mm^T_T)$$
such that $\Mm^T_T=\Mm^Z_Z|_T$. Note that we do not guarantee $(T,B_T,\Mm^T)/U$ is a g-sub-pair at this point.
\end{lem}
\begin{proof}
For any prime divisor $D$ on $Z$, we define $t_D:=t_D(X,B,\Mm;f)$. For any prime divisor $D$ on $X$, we define $D^v$ and $D^h$ to be the vertical$/Z$ part and the horizontal$/Z$ part of $D$ respectively.

Let $$\Theta:=B+\sum_{D\text{ is a component of }\Delta_Z}t_Df^*D.$$
By our construction, $\Delta\geq \Theta^{>0}$, and $(X,\Theta^{>0},\Mm)$ is log smooth and gdlt. Let $\Lambda$ be the set of all prime divisors $E$ over $X$, such that $f(E)\subset\Delta_Z$ and $\mult_E\Theta<1$. Then we may pick a real number $0<\epsilon\ll 1$ such that $(X,\Xi,\Mm):=(X,\Theta^{>0}+\epsilon\sum_{E\in\Lambda}E,\Mm)$ is glc. Now we may write
$$K_X+\Xi+\Mm_X\sim_{\Qq,Z}D:=D^h+D^v$$
where $D^h:=-B^{h,<0}\geq 0$ is horizontal$/Z$ and $D^v:=-\Theta^{v,<0}+\epsilon\sum_{E\in\Lambda}E$ is vertical$/Z$. By our construction, $D^v$ is very exceptional$/Z$. By assumption (2),
$$D^h|_F=-B^{h,<0}|_F=-B^{<0}|_F=N_{\sigma}((K_X+B^{>0}+\Mm_X)|_F).$$
By \cite[Lemma 2.19]{Hu20}, we may run a $(K_X+\Xi+\Mm_X)$-MMP$/Z$ which terminates with a good minimal model $(X',\Xi',\Mm)/Z$ of $(X,\Xi,\Mm)/Z$, such that the induced birational map $\phi: X\dashrightarrow X'$ contracts exactly $\Supp D$. 

For any prime divisor $D'$ on $X'$, we define $D'^v$ and $D'^h$ to be the vertical part$/Z$ and the horizontal$/Z$ part of $D'$ respectively. Let $B'$ be the strict transform of $B$ on $X'$ and $f': X'\rightarrow Z$ the induced contraction, then $$\Xi'=B'+\sum_{D\text{ is a component of }\Delta_Z}t_Df'^*D.$$
Since $f$ satisfies Condition \ref{cond: weak-semistabilized}, we have
$$\Xi'^v=(B'+\sum_{D\text{ is a component of }\Delta_Z}t_Df'^*D)^v=\sum_{D\text{ is a component of }\Delta_Z}f'^*D=f'^*\Delta_Z.$$
Since $\mult_SB=1$, $S\not\subset\Supp D$, hence $S$ is not contracted by $\phi$. We may let $S'$ be the strict transform of $S$ on $X'$ and $\tilde S'$ the normalization of $S'$. 

Let $f_{\tilde S'}:=f'|_{\tilde S'}: \tilde S'\rightarrow T$ be the induce morphism and $K_{\tilde S'}+\Xi_{\tilde S'}+\Mm^{S}_{\tilde S'}:=(K_{X'}+\Xi'+\Mm_{X'})|_{\tilde S'}$ the adjunction formulas. We let $\Delta_T:=\Delta_Z$ if $T=Z$, and $\Delta_T:=(\Delta_Z-T)|_T$ if $T$ is a prime divisor on $Z$. Then we have
\begin{equation}\label{equ: 1}
\Xi_{\tilde S'}\geq (B'+\sum_{D\text{ is a component of }\Delta_Z} t_Df'^*D-S')|_{\tilde S'}\geq \sum_{D_T\text{ is a component of }\Delta_T}f_{\tilde S'}^*D_T=f_{\tilde S'}^*\Delta_T.
\end{equation}
For any prime divisor $D_T$ on $T$, we define $t'_{D_T}:=t_{D_T}(S,B_S,\Mm^S;f|_S)$. By Condition \ref{cond: weak-semistabilized}(7), $f|_S$ has equi-dimensional and reduced fibers. By Lemma \ref{lem: hu20 3.33 gpair}, the discriminant part $\Supp B'_T\subset\Delta_T$, hence $$B'_T=\sum_{D_T\text{ is a component of }\Delta_T}(1-t'_{D_T})D_T.$$
Let $B_T:=B_Z$ if $T=Z$ and $B_T:=(B_Z-T)|_T$ if $T$ is a prime divisor on $Z$. For any component $D_T$ of $\Delta_T$, we may associate a component $D:=r(D_T)$ of $\Delta_Z$ to $\Delta_T$, such that $D|_T=D_T$, and we have $\mult_{D_T}B_T=\mult_DB_Z=1-t_D$. In particular, we have
$$B_T=\sum_{D_T\text{ is a component of }\Delta_T}(1-t_{r(D_T)})D_T.$$
Since $$(X,B+\sum_{D\text{ is a component of }\Delta_Z}t_Df^*D,\Mm)$$ is sub-glc, 
$$(S,B_S+\sum_{D_T\text{ is a component of }\Delta_T}t_{r(D_T)}(f|_S)^*D_T,\Mm^S)$$
is sub-glc. Thus $t_{r(D_T)}\leq t'_{D_T}$, which implies that $B_T\geq B'_T$. On the other hand, let
$$K_{\tilde S'}+B_{\tilde S'}+\Mm^S_{\tilde S'}:=(K_{X'}+B'+\Mm_{X'})|_{\tilde S'},$$
then $$\Xi_{\tilde S'}=B_{\tilde S'}+\sum_{D_T\text{ is a component of }\Delta_T}t_{r(D_T)}f_{\tilde S'}^*D_T.$$
By \ref{equ: 1}, we have
$$B_{\tilde S'}\geq\sum_{D_T\text{ is a component of }\Delta_T}(1-t_{r(D_T)})f_{\tilde S'}^*D_T,$$
hence
$$t'_{D_T}=t_{D_T}(\tilde S',B_{\tilde S'},\Mm^S;f_{\tilde S'})\leq t_{r(D_T)}$$
for any component $D_T$ of $\Delta_T$. Thus $B'_T\geq B_T$, hence $B_T=B'_T$, and the lemma follows.
\end{proof}

The following lemma is similar to \cite[Proposition 3.30]{Hu20}:
\begin{lem}\label{lem: hu20 3.30}
Let $(X,B,\Mm)/U$ be a glc $\Qq$-g-pair and $f: (X,B,\Mm)\rightarrow Z$ a good glc-trivial morphism$/U$. Then there exists a birational morphism $h_Z: Z'\rightarrow Z$ and a $\Qq$-factorial glc $\Qq$-g-pair $(X',B',\Mm)$ associated with a B-birational contraction $h: (X',B',\Mm)\dashrightarrow (X,B,\Mm)$, such that the induced map $f': (X',B',\Mm)\dashrightarrow Z'$ is a good gdlt-trivial morphism$/U$. Moreover, if $f$ is a contraction, then $f'$ is a contraction.
\end{lem}
\begin{proof}
Let $(Z,B_Z,\Mm^Z)$ be the $\Qq$-g-pair induced by $f: (X,B,\Mm)\rightarrow Z$, $h_Z: Z'\rightarrow Z$ be a gdlt modification of $(Z,B_Z,\Mm^Z)$, and $K_{Z'}+B_{Z'}+\Mm^Z_{Z'}:=h^*_Z(K_Z+B_Z+\Mm^Z_Z)$.

Let $p: W\rightarrow X$ be a log resolution of $(X,\Supp B)$ such that $\Mm$ descends to $W$ and the induced map $f_W: W\dashrightarrow Z'$ is a morphism. Let $B_W:=p^{-1}_*B+\Supp\Exc(p)$, then $(W,B_W)$ is log smooth and $(W,B_W,\Mm)$ is $\Qq$-factorial gdlt.

Let $X\xrightarrow{\tau}\tilde Z\xrightarrow{\gamma} Z$ and $W\xrightarrow{\tau'}\tilde Z'\xrightarrow{\gamma'} Z'$ be Stein factorizations. Since $h_Z$ is a birational morphism, the induced map $h_{\tilde Z}: \tilde Z'\rightarrow\tilde Z$ is a birational morphism. We let $(\tilde Z,B_{\tilde Z},\Mm^{\tilde Z})/U$ and $(\tilde Z',B_{\tilde Z'},\Mm^{\tilde Z})/U$ be the sub-glc pairs such that $K_{\tilde Z}+B_{\tilde Z}+\Mm^{\tilde Z}_{\tilde Z}:=\gamma^*(K_Z+B_Z+\Mm^Z_Z)$ and $K_{\tilde Z'}+B_{\tilde Z'}+\Mm^{\tilde Z}_{\tilde Z'}:=\gamma'^*(K_{Z'}+B_{Z'}+\Mm^Z_{Z'})$ where $\gamma^*\Mm^Z=\Mm^{\tilde Z}$ and $\gamma'^*\Mm^{Z'}=\Mm^{\tilde Z'}$. In particular we have $\Mm^{\tilde Z}=\Mm^{\tilde Z'}$ since $\Mm^{Z}=\Mm^{Z'}$. Notice that $(\tilde Z,B_{\tilde Z},\Mm^{\tilde Z})/U$ is exactly the glc pair induced by the canonical bundle formula of $\tau$ since $f$ is good. Therefore $(\tilde Z',B_{\tilde Z'},\Mm^{\tilde Z})/U$ is also the sub-glc pair induced by the canonical formula of $\tau'$, hence $f'$ is good. Possibly replacing $W$ with a higher model, for any prime divisor $D$ on $\tilde Z'$ that is exceptional$/\tilde Z$, we have that
\begin{itemize}
    \item $a(D,\tilde Z,B_{\tilde Z},\Mm^{\tilde Z})=0$, and
    \item $D$ is dominated by a component $S\subset\lfloor B_W\rfloor$ such that $a(S,X,B,\Mm)=0$
\end{itemize}
(cf. \cite[Version 2, Theorem 5.1]{HL21}). Let $E:=(K_W+B_W+\Mm_W)-p^*(K_X+B+\Mm_X)$ and let $E^h$ and $E^v$ be the horizontal$/\tilde Z'$ part and the vertical $/\tilde Z'$ part of $E$ respectively. Let $F$ be a general fiber of $\tau'$. Since $E^h$ is exceptional$/X$, $E^h=N_{\sigma}((K_W+B_W+\Mm_W)|_F)$.
\begin{claim}\
$E^v$ is very exceptional$/\tilde Z'$.
\end{claim}
\begin{proof}
Since $E^v$ is exceptional$/X$, $E^v$ is very exceptional$/\tilde Z$. Let $E^v=\sum_{i}a_iE_i$, then for any irreducible closed subset $T\subset f(\Supp E^v)$, we define $E_T:=\sum_{f(E_i)=T}a_iE_i$. There are three cases:
\begin{enumerate}
    \item $T$ is a codimension $\geq 2$ subset of $\tilde Z$. In this case, $E_T$ is very exceptional$/\tilde Z$.
    \item $T$ and $h_{\tilde Z}(T)$ are divisors. In this case, since $E_T$ is very exceptional$/\tilde Z$, $E_T$ is also very exceptional$/\tilde Z'$.
    \item $T$ is a $h_{\tilde Z}$-exceptional divisor. In this case, $T$ is dominated by a component $S\subset\lfloor B_W\rfloor$ such that $a(S,X,B,\Mm)=0$. Since $S\not\subset\Supp E^v$, $E_T$ is very exceptional$/\tilde Z'$.
\end{enumerate}
Thus $E^v=\sum_{T}E_T$ is  very exceptional$/\tilde Z'$ by definition.
\end{proof}
\noindent\textit{Proof of Lemma \ref{lem: hu20 3.30} continued}. By \cite[Lemma 2.19]{Hu20}, we may run a $(K_W+B_W+\Mm_W)$-MMP$/\tilde Z'$ which terminates with a good minimal model $(X',B',\Mm)/\tilde Z'$ such that the induced birational map $W\dashrightarrow X'$ contracts exactly $E$. Therefore, the induced birational map $h: (X',B',\Mm)\dashrightarrow (X,B,\Mm)$ is a B-birational contraction, and the induced morphism $f': (X',B',\Mm)\rightarrow Z'$ is a gdlt-trivial morphism$/U$. 

If $f$ is a contraction, then $Z=\tilde Z$ and $Z'=\tilde Z'$, hence $f'$ is a contraction.
\end{proof}

The statement and the proof of the following theorem is similar to \cite[Theorem 3.38]{Hu20}.

\begin{thm}\label{thm: hu20 3.38 gpair}
%Let $(X,B,\Mm)/U$ be a gdlt $\Qq$-g-pair and 
Let $(X,B,\Mm)/U$ be a gdlt $\Qq$-g-pair, $f: (X,B,\Mm)\rightarrow Z$ a good gdlt-trivial morphism$/U$, $S$ a component of $\lfloor B\rfloor$, and $T=Z$ or a component of $\lfloor B_Z\rfloor$, such that $f(S)=T$. 

Let  $(S,B_S,\Mm^S)/U$ be the gdlt $\Qq$-g-pair induced by the adjunction $K_S+B_S+\Mm^S_S:=(K_X+B+\Mm_X)|_S$, $(T,B_T,\Mm^T)/U$ the glc $\Qq$-g-pair induced by the glc-trivial morphism$/U$ $f|_S: (S,B_S,\Mm^S)\rightarrow T$, and $(T,B'_T,\NN^T)/U$ the gdlt $\Qq$-g-pair induced by the adjunction $K_T+B'_T+\NN^T_T:=(K_Z+B_Z+\Mm^Z_Z)|_T$. Then:
\begin{enumerate}
\item $f|_S: (S,B_S,\Mm^S)\rightarrow T$ is good, and
    \item  $B_T'=B_T$, hence we may choose $\Mm^T$ and $\Mm^Z$ such that $\NN^T=\Mm^T$. In particular, $f|_S: (S,B_S,\Mm^S)\rightarrow T$ is a gdlt-trivial morphism$/U$.
\end{enumerate}
\end{thm}
\begin{proof}
%When $\dim X=1$ there is nothing to prove. When $\dim X\geq 2$, by induction on dimensions, we may assume that $S$ is a prime divisor on $X$ and either $T=Z$ or $T$ is a prime divisor on $Z$.  % By Lemma \ref{lem: hu20 3.30} we may assume that $f: (X,B,\Mm)\rightarrow Z$ is a gdlt-trivial morphism$/U$.
\noindent\textbf{Step 1}. In this step we reduce to the case when $f$ is a contraction. Let $X\xrightarrow{\tau}\tilde Z\xrightarrow{\gamma} Z$ be the Stein factorization of $f$ and let $S\xrightarrow{\tau_S}\tilde T\xrightarrow{\gamma_S}T$ be the induced factorization, where $\tilde T$ is the normalization of the image of $S$ on $\tilde Z$. Since $\gamma_S$ is finite and $f|_S: (S,B_S,\Mm^S)\rightarrow T$ is a glc-trivial morphism$/U$, $\tau_S: (S,B_S,\Mm^S)\rightarrow\tilde T$ is a glc-trivial morphism$/U$. Let $(\tilde T,B_{\tilde T},\Mm^{\tilde T})/U$ be the glc g-pair induced by $\tau_S: (S,B_S,\Mm^S)\rightarrow\tilde T$, then $\Mm^T=\frac{1}{\deg\gamma_S}(\gamma_S)_*\Mm^{\tilde T}$. Let $(\tilde Z,B_{\tilde Z},\Mm^{\tilde Z})/U$ be the g-pair induced by a canonical bundle formula
$$K_X+B+\Mm_X\sim_{\Qq}\tau^*(K_{\tilde Z}+B_{\tilde Z}+\Mm^{\tilde Z}_{\tilde Z})$$
such that $K_{\tilde Z}+B_{\tilde Z}+\Mm^{\tilde Z}_{\tilde Z}=\gamma^*(K_Z+B_Z+\Mm^Z_Z)$, then $\Mm^{\tilde Z}=\gamma_S^*\Mm^Z$ since $f$ is good. Let $(\tilde T,B_{\tilde T}',\NN^{\tilde T})/U$ be the glc g-pair induced by the adjunction $K_{\tilde T}+B_{\tilde T}'+\NN^{\tilde T}_{\tilde T}:=(K_{\tilde Z}+B_{\tilde Z}+\Mm^{\tilde Z}_{\tilde Z})|_{\tilde T}$.

Notice that $\NN^{\tilde T}=\Mm^{\tilde Z}|_{\tilde T}=\gamma_S^*\Mm^Z|_T=\gamma_S^*\NN^T$. Then by \cite[Lemma 3.22]{Hu20}, we only need to show that $\tau_S: (S,B_S,\Mm^S)\rightarrow\tilde T$ is good and $\NN^{\tilde T}=\Mm^{\tilde T}$. By Lemma \ref{lem: hu20 3.30}, there exists a gdlt modification $\bar Z\rightarrow\tilde Z$ of $(\tilde Z,B_{\tilde Z},\Mm^{\tilde Z})$ and a gdlt-trivial morphism$/U$ $\bar f: (\bar X,\bar B,\Mm)\rightarrow\bar Z$ such that the induced birational map $h: (\bar X,\bar B,\Mm)\dashrightarrow (X,B,\Mm)$ is B-birational and $\bar f$ is a contraction. Possibly replacing $f: (X,B,\Mm)\rightarrow Z$ with $\bar f: (\bar X,\bar B,\Mm)\rightarrow\bar Z$, we may assume that $f$ is a contraction, and $f: (X,B,\Mm)\rightarrow Z$ is a glc-trivial fibration$/U$.

\medskip

\noindent\textbf{Step 2}. In this step we reduce to the case when $f|_S$ is a contraction. Let $S\xrightarrow{\tau_T}T'\xrightarrow{\gamma_T}T$ be the Stein factorization of $f|_S$. Either $T=Z$, then we choose $Z'$ to be $T'$. Or $T$ is a prime divisor in $Z$, then by \cite[Theorem A.1]{Hu20}, there exists a finite morphism $\gamma_Z: Z'\rightarrow Z$ together with a prime divisor $\hat{T}$ on $Z'$ such that the induced morphism $\nu: T'\rightarrow\hat{T}$ is the normalization. Then there exists a an induced finite morphism $\gamma: X'\rightarrow X$ and an induced glc-trivial fibration$/U$ $f': (X',B',\Mm')\rightarrow Z'$. In either case, by \cite[Fact ($\clubsuit$)]{Hu20}, there exists a component $S'$ of $S\times_TT'$ which is isomorphic to $S$. Let $(S',B_{S'},\Mm^{S'})/U$ be the g-pair induced by $K_{S'}+B_{S'}+\Mm^{S'}_{S'}:=(K_{X'}+B'+\Mm'_{X'})|_{S'}$. Then we have an induced contraction $f'_{S'}:=f'|_{S'}: S'\rightarrow T'$, such that $f'_{S'}: (S',B_{S'},\Mm^{S'})\rightarrow T'$ is a glc-trivial fibration$/U$. We let 
$$K_{S'}+B_{S'}+\Mm^{S'}_{S'}\sim_{\Qq}f'^*_{S'}(K_{T'}+B_{T'}+\Mm^{T'}_{T'})$$
and
$$K_{X'}+B'+\Mm'_{X'}\sim_{\Qq}f'^*(K_{Z'}+B_{Z'}+\Mm^{Z'}_{Z'})$$
be canonical bundle formulas, and let 
$$K_{T'}+B'_{T'}+\NN^{T'}_{T'}:=(K_{Z'}+B_{Z'}+\Mm^{Z'}_{Z'})|_{T'}$$
be the adjunction formulas. Since the induced g-pair structure of glc-trivial fibrations are compatible with generically finite base change \cite[Chapter 6, Remark 7]{Fil19} (see also \cite[Remark 4.3]{Fil20}), we can pick canonical bundle formulas such that $\Mm^{Z'}=\gamma_Z^*\Mm^Z$ and $\NN^{T}=\frac{1}{\deg\gamma_T}(\gamma_T)_*\NN^{T'}$. By \cite[Lemma 3.22]{Hu20}, we only need to prove $\Mm^{T'}=\NN^{T'}$ which will imply that $f|_S: (S,B_S,\Mm^S)\rightarrow T$ is good. By Lemma \ref{lem: hu20 3.30}, there exists a gdlt modification $Z''\rightarrow Z'$ of $(Z',B_{Z'},\Mm^{Z'})$ and a gdlt-trivial fibration$/U$ $f'':(X'',B'',\Mm')\rightarrow Z''$ such that the induced birational map $(X'',B'',\Mm')\dashrightarrow (X',B',\Mm')$ is a B-birational contraction. Possibly replacing $f: (X,B,\Mm)\rightarrow Z$ with $f'': (X'',B'',\Mm')\rightarrow Z''$, and replacing $S,T$ with the strict transforms of $S'$ and $T'$ on $X''$ and $Z''$ respectively, we may assume that $f|_S$ is a contraction.

\medskip

\noindent\textbf{Step 3}. We finish the proof in this step. By \cite[Thereoms B.6, B.7]{Hu20}, there exist birational morphisms $\pi: X_1\rightarrow X$ and $\phi: Z_1\rightarrow Z$ such that the induced map $f_1: X_1\rightarrow Z_1$ is a flat morphism, a Galois finite morphism $\psi: Z_2\rightarrow Z_1$ with induced morphisms $\mu: X_2\rightarrow X_1$ and $f_2: X_2\rightarrow Z_2$ where $X_2$ is the main component of $X_1\times_{Z_1}Z_2$, and two reduced divisors $\Delta$ on $X_2$ and $\Delta_Z$ on $Z_2$, such that $(X_2,B_2,\Mm^2)/U,(Z_2,B_{Z_2},\Mm^{2,Z})/U$, $f_2: X_2\rightarrow Z_2$, and $\Delta$ and $\Delta_Z$ satisfy Condition \ref{cond: weak-semistabilized}, where $K_{X_2}+B_2+\Mm^2_{X_2}:=(\pi\circ\mu)^*(K_X+B+\Mm_X)$. 

Let $S_2$ be the strict transform of $S$ on $X_2$, $T_2$ the strict transform of $T$ on $Z_2$, and $F$ a general fiber of $f_2$. Since the induced g-sub-pair structure of glc-trivial fibrations are compatible with generically finite base change \cite[Chapter 6, Remark 7]{Fil19} (see also \cite[Remark 4.3]{Fil20}), the theorem follows by applying Lemma \ref{lem: hu20 3.35 gpair} to $(X_2,B_2,\Mm^2)/U,(Z_2,B_{Z_2},\Mm^{2,Z})/U,f_2,\Delta,\Delta_Z$, $S_2$ and $F$.
\end{proof}

\begin{thm}\label{thm: b log abundance good case}
Let $(X,B,\Mm)/U$ be a glc $\Qq$-g-pair, $f: (X,B,\Mm)\rightarrow Z$ a good glc-trivial morphism$/U$, and $(Z,B_Z,\Mm^Z)/U$ a glc $\Qq$-g-pair induced by $f: (X,B,\Mm)\rightarrow Z$. Assume that $\Mm$ is $\bb$-log abundant$/U$ with respect to $(X,B,\Mm)$. Then $\Mm^Z$ is $\bb$-log abundant$/U$ with respect to $(Z,B_Z,\Mm^Z)$.
\end{thm}
\begin{proof}
By Lemma \ref{lem: hu20 3.30}, we may assume that $f: (X,B,\Mm)\rightarrow Z$ is a gdlt-trivial morphism$/U$. Since any component of $\lfloor B_Z\rfloor$ is dominated by a component of $\lfloor B\rfloor$ (cf. \cite[Version 2, Theorem 5.1]{HL21}), by Theorem \ref{thm: hu20 3.38 gpair} and induction on dimensions, we only need to prove that $\Mm^Z$ is $\bb$-abundant$/U$. By Theorem \ref{thm: hu20 3.38 gpair} again, we may assume that $(X,B,\Mm)$ does not have any glc center that is horizontal$/Z$. The theorem follows from Lemma \ref{lem: abundant gcbf klt ssample case}.
\end{proof}

\begin{lem}\label{lem: contraction log abundance gcbf}
Theorem \ref{thm: log abundance gcbf} holds when $f$ is a contraction.
\end{lem}
\begin{proof}
The $\Qq$-coefficient follows from Theorem \ref{thm: b log abundance good case} as any contraction is good. Thanks to Lemma \ref{lem: abundant is open in minimal Q affine V} and Theorem \ref{thm: cbf for R-coefficient}, the $\Rr$-coefficient case follows from the $\Qq$-coefficient case and the standard theory of Shokurov-type rational polytopes (cf. \cite[Proposition 3.16]{HL18}) and uniform rational polytopes (cf. \cite[Lemma 5.3]{HLS19}, \cite[Therem 1.4]{Che20}). See also \cite[Lemma 4.1]{HL19}.
\end{proof}

\begin{proof}[Proof of Theorem \ref{thm: log abundance gcbf}]
By Lemma \ref{lem: contraction log abundance gcbf} and taking the Stein factorization of $f$, we can assume that $f$ is a finite morphism. After birational base change we can also assume that $\Mm$ descends to $X$, $\Mm^Z$ descends to $Z$ and $Z$ is smooth.
Notice that $\Mm^Z=\frac{1}{\deg f}f_*\Mm$ is already nef and abundant by \cite[Lemma 2.12]{Hu20}, so we only need to prove that $\Mm_Z^Z|_T$ is abundant for any g-sub-lc center $T$ of $(Z,B_Z,\Mm^Z)/U$. Possibly by replacing $Z$ with higher models we can assume $T$ is a prime divisor on $Z$. Let $S_1,\dots,S_n$ be the components of $f^*T$, $r_i$ the ramification index of $f$ along $S_i$, $d_i$ be the degree of the finite morphism $f|_{S_i}:S_i\to T$, and $a_i:=\mult_{S_i}B_i$. Since $(X,B,\Mm)/U$ is a sub-glc g-sub-pair, $a_i\leq 1$ for each $i$.

By the Hurwitz formula, We have $K_X=f^*K_Z+\sum_i(r_i-1)S_i$ over a neighborhood of the generic point of $T$. By Definition-Lemma \ref{deflem: cons g-pair finite morphism}, over a neighborhood of the generic point of $T$, we have
$$T=B_Z=\frac{1}{\deg f}f_*(\sum_i (r_i+a_i-1)S_i)=\sum_i(r_i+a_i-1)d_iT.$$
Since $\sum_ir_id_i=\deg f$, $a_i=1$ for each $i$.

Therefore, each $S_i$ is a glc center of $(X,B,\Mm)$. Since $\Mm_X$ is $\bb$-abundant$/U$, $\Mm_X|_{S_i}$ is $\bb$-abundant over $U$. By the projection formula, $\Mm^Z_Z|_T=\sum_i\frac{1}{d_i}f_{S_i,*}\Mm_X|_{S_i}$. By \cite[Lemmas 2.12]{Hu20} and Lemma \ref{lem: equiv def for nef and abundant}(3), $\Mm^Z_Z|_T$ is abundant$/U$, which concludes the proof.
\end{proof}

\begin{proof}[Proof of Theorem \ref{thm: sub-adjunction log abundant}]
The proof is similar to the proof of \cite[Theorem 5.1]{HL19}. Let $V$ be a glc center of $(X,B,\Mm)$ with normalization $W$. Let $f: Y\rightarrow X$ be a gdlt modification of $(X,B,\Mm)$ such that there exists a prime divisor $S\subset B_Y\rfloor$ such that $f(S)=V$, where $K_Y+B_Y+\Mm_Y:=f^*(K_X+B+\Mm_X)$. Let $W_Y$ be a glc center of $(Y,B_Y,\Mm)$ which is minimal with respect to inclusion under the condition $f(W_Y)=V$. Since $(Y,B_Y,\Mm)$ is gdlt, by repeatedly applying adjunction, we get a gdlt g-pair $(W_Y,B_{W_Y},\Mm^{W_Y})/U$ such that
$$K_{W_Y}+B_{W_Y}+\Mm^{W_Y}_{W_Y}:=(K_Y+B_Y+\Mm_Y)_{W_Y}.$$
Since $\Mm$ is $\bb$-log abundant$/U$ with respect to $(X,B,\Mm)$, $\Mm^{W_Y}$ is $\bb$-log abundant$/U$ with respect to $(W_Y,B_{W_Y},\Mm^{W_Y})$. By construction, there exists a naturally induced projective surjective morphism $f_W: W_Y\rightarrow W$ such that $K_{W_Y}+B_{W_Y}+\Mm^{W_Y}_{W_Y}\sim_{\Rr,W}0$. By Theorem \ref{thm: log abundance gcbf}, there exists a glc g-pair $(W,B_W,\Mm^W)/U$ induced by $f_W: (W_Y,B_{W_Y},\Mm^{W_Y})\rightarrow W$, such that $\Mm^W$ is $\bb$-log abundant$/U$ with respect to $(W,B_W,\Mm^W)$. $(W,B_W,\Mm^W)/U$ is the glc g-pair we desire.
\end{proof}

\section{Non-vanishing and abundance}

In this section, we prove Theorem \ref{thm: klt abundant to glc good nef abundant} and Corollary \ref{cor: log abundant m dim 3 abundance}.

\begin{proof}[Proof of Theorem \ref{thm: klt abundant to glc good nef abundant}] It is clear that Theorem \ref{thm: klt abundant to glc good nef abundant} holds when $d=1$. If $d\geq 2$, by induction, we may further assume that Theorem \ref{thm: klt abundant to glc good nef abundant} holds in dimension $\leq d-1$. Possibly replacing $(X,B,\Mm)$ with a gdlt modification, we may assume that $(X,B,\Mm)$ is $\Qq$-factorial gdlt. In particular, if we are in case (2), then we only need to prove that $K_X+B+\Mm_X$ is abundant$/U$ and the log abundance$/U$ will follow from induction on dimension by applying adjunction to glc centers.

We may assume that $K_X+B+\Mm_X$ is pseudo-effective$/U$, otherwise there is nothing to prove. 

If $K_X+B-\epsilon_0\lfloor B\rfloor+\Mm_X$ is pseudo-effective$/U$ for some positive real number $\epsilon_0$, then since $\Mm$ is $\bb$-abundant$/U$, we may pick $0<\epsilon\ll\epsilon_0$ and pick $0\leq \Delta\sim_{\Rr,U}B-\epsilon\lfloor B\rfloor+\Mm_X$ such that $(X,\Delta)$ is klt. By the non-vanishing conjecture for klt pairs in dimension $d$, $K_X+\Delta\sim_{\Rr,U}D\geq 0$ for some $\Rr$-divisor $D$. Thus $K_X+B+\Mm_X\sim_{\Rr,U}D+\epsilon\lfloor B\rfloor\geq 0$. Moreover, if we are in case (2), then $K_X+B-\epsilon\lfloor B\rfloor+\Mm_X$ is abundant$/U$, and $K_X+B+\Mm_X$ is abundant$/U$ by Lemma \ref{lem: limit of nef and abundant divisors is nef and abundant}. Therefore, we may assume that $K_X+B-\epsilon\lfloor B\rfloor+\Mm_X$ is not pseudo-effective$/U$ for any positive real number $\epsilon$. 

For any real number $0<\epsilon\ll 1$, since $(X,B-\epsilon\lfloor B\rfloor,\Mm)$ is gklt, by \cite[Lemma 4.4]{BZ16}, we may run a $(K_X+B-\epsilon\lfloor B\rfloor+\Mm_X)$-MMP$/U$ which terminates with a Mori fiber space$/U$: $f_{\epsilon}: X_{\epsilon}\rightarrow Z_{\epsilon}$. Let $B_{\epsilon}$ be the birational transform of $B$ on $X_{\epsilon}$ for each $0<\epsilon\ll 1$. Since $K_X+B+\Mm_X$ is pseudo-effective$/U$, $K_{X_{\epsilon}}+B_{\epsilon}+\Mm_{X_{\epsilon}}$ is pseudo-effective$/U$ for any $0<\epsilon\ll 1$. Since $\rho(X_{\epsilon}/Z_{\epsilon})=1$, for any $0<\epsilon\ll 1$, there exists $\delta_{\epsilon}\in [0,\epsilon)$ such that $K_{X_{\epsilon}}+B_{\epsilon}-\delta_{\epsilon}\lfloor B_{\epsilon}\rfloor+\Mm_{X_{\epsilon}}\sim_{\Rr,Z_{\epsilon}}0$. By \cite[Theorem 1.5]{BZ16}, for any $0<\epsilon\ll 1$, $(X_{\epsilon},B_{\epsilon},\Mm)$ is glc, hence $(X_{\epsilon},B_{\epsilon}-\delta_{\epsilon}\lfloor B_{\epsilon}\rfloor,\Mm)$  is glc. By \cite[Theorem 1.6]{BZ16}, for any $0<\epsilon\ll 1$, $\delta_{\epsilon}=0$ and $K_{X_{\epsilon}}+B_{\epsilon}+\Mm_{X_{\epsilon}}\sim_{\Rr,Z_{\epsilon}}0$. In the following, we will fix $0<\epsilon\ll 1$ and denote $X':=X_{\epsilon}$, $Z':=Z_{\epsilon}$, and $B_{\epsilon}:=B'$.

Let $p: \bar X\rightarrow X$ and $q: \bar X\rightarrow X'$ be a common resolution such that $\Mm$ descends to $\bar X$, $p$ is a log resolution of $(X,\Supp B)$, and $q$ is a log resolution of $(X',\Supp B')$. Then we may write
$$K_{\bar X}+\bar B+\Mm_{\bar X}=p^*(K_X+B+\Mm_X)+E$$
for some $p$-exceptional $\Rr$-divisor $E\geq 0$ and gdlt g-pair $(\bar X,\bar B,\Mm)/U$ such that $\bar B\wedge E=0$, and write
$$K_{\bar X}+\bar B+\Mm_{\bar X}=q^*(K_{X'}+B'+\Mm_{X'})+F-G$$
for some $q$-exceptional $\Rr$-divisors $F\geq 0,G\geq 0$ such that $F\wedge G=0$. By \cite[Proposition 3.8]{HL18}, we may run a $(K_{\bar X}+\bar B+\Mm_{\bar X})$-MMP$/X'$ which terminates with a log minimal model $(\bar X',\bar B',\Mm)/X'$ of $(\bar X,\bar B,\Mm)$, such that $$K_{\bar X'}+\bar B'+\Mm_{\bar X'}+G_{\bar X'}=q'^*(K_{X'}+B'+\Mm_{X'})\sim_{\Rr,Z'}0,$$ where $G_{\bar X'}$ is the strict transform of $G$ on $\bar X'$ and $q':\bar X'\rightarrow X'$ is the induced morphism.

Since $K_X+B+\Mm_X$ is pseudo-effective$/U$, $K_{\bar X'}+\bar B'+\Mm_{\bar X'}$ is pseudo-effective$/U$ by our construction. Since $\rho(X'/Z')=1$, $X'$ is $\Qq$-factorial klt, and $K_{X'}+B'+\Mm_{X'}\sim_{\Rr,Z'}0$, we may find $0\leq\Delta\sim_{\Rr,Z'}B'+\Mm_{X'}$ such that $(X,\Delta)$ is klt. By \cite[Proposition 3.2]{Has20}, there exists a $(K_{\bar X'}+\bar B'+\Mm_{\bar X'})$-negative birational contraction$/Z'$: $\bar X'\dashrightarrow X''$ such that $K_{X''}+B''+\Mm_{X''}$ is semi-ample$/Z'$, where $B''$ is the strict transform of $\bar B'$ on $X''$. In particular, $\Mm$ is $\bb$-log abundant$/U$ with respect to $(X'',B'',\Mm)$. We let $f: X''\rightarrow Z$ be the contraction$/Z'$ defined by $K_{X''}+B''+\Mm_{X''}$. Since $K_{X'}+B'+\Mm_{X'}\sim_{\Rr,Z'}0$, the restriction of $K_{X''}+B''+\Mm_{X''}$ to a general fiber of $X''\rightarrow Z'$ is $\Rr$-linearly equivalent to $0$. Thus $Z$ is birational to $Z'$. 

Since $K_{X''}+B''+\Mm_{X''}\sim_{\Rr,Z}0$ and $\Mm$ is $\bb$-log abundant$/U$ with respect to $(X'',B'',\Mm)$, by Theorem \ref{thm: log abundance gcbf}, we have a canonical bundle formula
$$K_{X''}+B''+\Mm_{X''}\sim_{\Rr}f^*(K_Z+B_Z+\Mm^Z_Z)$$
for some glc g-pair $(Z,B_Z,\Mm^Z)/U$ such that $\Mm^Z$ is $\bb$-log abundant$/U$ with respect to $(Z,B_Z,\Mm^Z)/U$. Since $\dim Z\leq d-1$, by Theorem \ref{thm: klt abundant to glc good nef abundant} in dimension $\leq d-1$, $|K_Z+B_Z+\Mm^Z_Z/U|\not=\emptyset$, and if we are in case (2), then $K_Z+B_Z+\Mm^Z_Z$ is abundant$/U$. Thus $|K_{X''}+B''+\Mm_{X''}/U|\not=\emptyset$, and if we are in case (2), then  $K_{X''}+B''+\Mm_{X''}$ is abundant$/U$. Theorem \ref{thm: klt abundant to glc good nef abundant} follows from our construction of $(X'',B'',\Mm)/U$.
\end{proof}

\begin{proof}[Proof of Corollary \ref{cor: log abundant m dim 3 abundance}]
It immediately follows from Theorem \ref{thm: klt abundant to glc good nef abundant} and the abundance theorem for klt pairs in dimension $\leq 3$ (cf. \cite{KMM94})
\end{proof}

\section{Example}

The following example shows that, given a canonical bundle formula
$K_X+B+\Mm_X\sim_{\Qq}f^*(K_Z+B_Z+\Mm^Z_Z)$
with $(X,B,\Mm)$ glc but not gklt, it is possible that $\Mm^Z$ is not $\bb$-abundant even if $\Mm=\overline{\Mm_X}$ is $\bb$-abundant.
\begin{ex}[=Example \ref{ex: intro example ab not log ab}]\label{ex: counter-example not b-log abundant}
Let $E$ be a smooth elliptic curve, $L$ an ample divisor on $E$, $V:=\mathcal{O}_E\oplus\mathcal{O}_E(L)$, $X:=\mathbb P_E(V)$, and $f: X\rightarrow E$ the induced contraction. Let $v$ be the class of a fiber of $f$ and let $h$ be the class of sheaf $\mathcal{O}_X(1)$ in $\Pic(X)$. Then there are two sections $E_1,E_2$ on $X$, such that $[E_1]=h-f^*(L)$ is the negative section and $[E_2]=h$ is the positive section. Notice that $c_2(V)=0$ and $c_1(V)=c_1(L)$, we have $h^2=f^*c_1(L)\cdot h=\deg_E(L)$. Since $L$ is ample, $E_2$ is nef and big as $E^2=h^2>0$. Since $\omega_{X/E}=f^*(\det V)\otimes\mathcal{O}_X(-2)$, we have $K_X+E_1+E_2\sim 0$, hence $K_X+E_1+E_2\sim f^*K_E$. 

Let $P$ be a non-torsion numerically trivial divisor on $E$. Then $(X,B:=E_1,\Mm:=\overline{E_2+f^*P})$ is a projective glc g-pair and $K_X+B+\Mm_X\sim_{\Qq}f^*(K_E+P)$. In particular, $K_X+B+\Mm_X\sim_{\Qq,E}0$. Since $(X,E_1)$ is log smooth over $E$, $(E,B_E:=0,\Mm^E:=\overline{P})$ is a projective g-pair induced by a canonical bundle formula. Since $E_2$ is big and nef, $E_2+f^*P$ is big and nef, hence abundant. Thus $\Mm$ is $\bb$-abundant. However, $P$ is not abundant, hence $\Mm^E$ is not $\bb$-abundant. Moreover, it is clear that $K_X+B+\Mm_X$ is pseudo-effective but not effective. In particular, $K_X+B+\Mm_X$ is not abundant.
\end{ex}


\begin{thebibliography}{99}

\bibitem[AK00]{AK00} D. Abramovich and K. Karu, \textit{Weak semistable reduction in characteristic 0}, Invent. Math. \textbf{139} (2000), no. 2, 241--273.

%\bibitem[AHK07]{AHK07} V. Alexeev, C. D. Hacon, and Y. Kawamata, \textit{Termination of (many) $4$-dimensional log flips}, Invent. Math. \textbf{168} (2007), no. 2, 433--448.

\bibitem[Amb99]{Amb99} F. Ambro, \textit{The Adjunction Conjecture and its applications}, arXiv:math/9903060v3.

%\bibitem[Amb03]{Amb03} F. Ambro, \textit{Quasi-log varieties}, Tr. Mat. Inst. Steklova \textbf{240} (2003), Biratsion. Geom. Linein. Sist. Konechno Porozhdennye Algebry, 220--239; translation in Proc. Steklov Inst. Math. 2003, no. 1 (240), 214--233.

\bibitem[Amb05]{Amb05} F. Ambro, \textit{The moduli b-divisor of an lc-trivial fibration}, Compos. Math. \textbf{141} (2005), no. 2, 385--403.

%\bibitem[Bir07]{Bir07} C. Birkar, \textit{Ascending chain condition for log canonical thresholds and termination of log flips}, Duke Math. J. \textbf{136} (2007), no. 1, 173--180. 

%\bibitem[Bir12a]{Bir12a} C. Birkar, \textit{Existence of log canonical flips and a special LMMP}, Pub. Math. IHES., \textbf{115} (2012), 325--368.

%\bibitem[Bir12b]{Bir12b} C. Birkar, \textit{On existence of log minimal models and weak Zariski decompositions}, Math. Ann. \textbf{354} (2012), no. 2, 787--799.

%\bibitem[Bir18]{Bir18} C. Birkar, \textit{Log Calabi-Yau fibrations}, arXiv: 1811.10709v2.
	
\bibitem[Bir19]{Bir19} C. Birkar, \textit{Anti-pluricanonical systems on Fano varieties}. Ann. of Math. (2), \textbf{190} (2019), 345--463.
	
%\bibitem[Bir20a]{Bir20a} C. Birkar, \textit{Geometry and moduli of polarised varieties}, arXiv: 2006.11238v1.

\bibitem[Bir20a]{Bir20a} C. Birkar, \textit{Generalised pairs in birational geometry}, arXiv: 2008.01008v2.

\bibitem[Bir20b]{Bir20b} C. Birkar, \textit{On connectedness of non-klt loci of singularities of pairs}, arXiv: 2010.08226v1.

\bibitem[Bir21a]{Bir21a} C. Birkar, \textit{Singularities of linear systems and boundedness of Fano varieties}, Ann. of Math. \textbf{193} (2021), no. 2, 347--405.


\bibitem[Bir21b]{Bir21b} C. Birkar, \textit{Boundedness and volume of generalised pairs}, arXiv: 2103.14935v2.

\bibitem[BCHM10]{BCHM10} C. Birkar, P. Cascini, C. D. Hacon and J. M\textsuperscript{c}Kernan, \textit{Existence of minimal models for varieties of log general type}, J. Amer. Math. Soc. \textbf{23} (2010), no. 2, 405--468.

%\bibitem[BDCS20]{BDCS20} C. Birkar, G. Di Cerbo, and R. Svaldi, \textit{Boundedness of elliptic Calabi-Yau varieties with a rational section}, arXiv: 2010.09769v1.

\bibitem[BH14]{BH14} C. Birkar and Z. Hu, \textit{Polarized pairs, log minimal models, and Zariski decompositions}, Nagoya Math. J. \textbf{215}(2014), 203--224.


\bibitem[BZ16]{BZ16} C. Birkar and D.-Q. Zhang, \textit{Effectivity of Iitaka fibrations and pluricanonical systems of polarized pairs}, Pub. Math. IHES., \textbf{123} (2016), 283--331.

\bibitem[Che20]{Che20} G. Chen, \textit{Boundedness of $n$-complements for generalized pairs}, arXiv: 2003.04237v2.

%\bibitem[CT20]{CT20} G. Chen and N. Tsakanikas, \textit{On the termination of flips for log canonical generalized pairs}, arXiv: 2011.02236v1.

%\bibitem[CX20]{CX20} G. Chen and Q. Xue, \textit{Boundedness of $(\epsilon,n)$-Complements for projective generalized pairs of Fano type}, arXiv: 2008.07121v1.

\bibitem[Cho08]{Cho08} R. Choi, \textit{The geography of log models and its applications}, PhD Thesis, Johns Hopkins University (2008).

%\bibitem[Fil20a]{Fil20a} S. Filipazzi, \textit{Boundedness of Log Canonical Surface Generalized Polarized Pairs}, Taiwanese J. Math. \textbf{22} (2018), no.4, 813--850.

\bibitem[Fil19]{Fil19} S. Filipazzi, \textit{Generalized pairs in birational geometry}, 2019. PhD thesis, University of Utah.

\bibitem[Fil20]{Fil20} S. Filipazzi, \textit{On a generalized canonical bundle formula and generalized adjunction}, Ann. Sc. Norm. Super. Pisa Cl. Sci. (5) Vol. XXI (2020), 1187--1221.

%\bibitem[Fil20]{Fil20} S. Filipazzi, \textit{On the boundedness of $n$-folds with $\kappa(X)=n-1$}, arXiv: 2005.05508v2. 

%\bibitem[FS20a]{FS20a} S. Filipazzi and R. Svaldi, \textit{Invariance of Plurigenera and boundedness for Generalized Pairs}, arXiv: 2005.04254v2.

\bibitem[FS20]{FS20} S. Filipazzi and R. Svaldi, \textit{On the connectedness principle and dual complexes for generalized pairs}, arXiv: 2010.08018v2.

%\bibitem[FW20]{FW20} S. Filipazzi and J. Waldron, \textit{Connectedness principle in characteristic $p>5$}, arXiv: 2010.08414v2.

\bibitem[FL19]{FL19} E. Floris and V. Laz\'ic, \textit{On the B-Semiampleness Conjecture}, \'Epijournal de G\'eom\'etrie Alg\'ebrique, Vol. 3 (2019), Article Nr. 12.



%\bibitem[Fuj04]{Fuj04} O. Fujino, \textit{Termination of $4$-fold canonical flips}, Publ. Res. Inst. Math. Sci, \textbf{40} (2004), no. 1, 231--237.

%\bibitem[Fuj05]{Fuj05} O. Fujino, \textit{Addendum to “Termination of $4$-fold canonical flips”}, Publ. Res. Inst. Math. Sci. \textbf{41} (2005), no. 1, 252--257.

%\bibitem[Fuj10]{Fuj10} O. Fujino, \textit{On Kawamata’s theorem}, Classification of Algebraic Varieties, EMS Ser. of Congr. Rep., Eur. Math. Soc., Z\"urich (2010), 305--315.

%\bibitem[Fuj11]{Fuj11} O. Fujino, \textit{Foundations of the minimal model program}, MSJ Memoirs, \textbf{35}. Mathematical Society of Japan, Tokyo (2017).

%\bibitem[Fuj12]{Fuj12} O. Fujino, \textit{Base point free theorems: saturation, B-divisors, and canonical bundle formula}, Algebra Number Theory \textbf{6} (2012), no. 4, 797--823.

%\bibitem[Fuj13]{Fuj13} O. Fujino, \textit{A transcendental approach to Koll\'ar’s injectivity theorem II},  J. Reine Angew. Math. \textbf{681} (2013), 149--174.

%\bibitem[Fuj18]{Fuj18} O. Fujino, \textit{Fundamental properties of basic slc-trivial fibrations I}, to appear in Publ. Res. Inst. Math. Sci., arXiv: 1804.11134v3.


%\bibitem[Fuj19]{Fuj19} O. Fujino, \textit{Corrigendum: On subadditivity of the logarithmic Kodaira dimension}, arXiv: 1904.11639v3.

%\bibitem[Fuj21]{Fuj21} O. Fujino, \textit{Cone theorem and Mori hyperbolicity},  arXiv:2102.11986v1.

\bibitem[FG12]{FG12} O. Fujino and Y. Gongyo \textit{On canonical bundle formulas and subadjunctions}, Michigan Math. J. \textbf{61} (2012), no. 2, 255--264.

\bibitem[FG14]{FG14} O. Fujino and Y. Gongyo, \textit{On the moduli b-divisors of lc-trivial fibrations}, Ann. Inst. Fourier (Grenoble), \textbf{64} (2014), no. 4, 1721--1735.

%\bibitem[FG14]{FG14} O. Fujino and Y. Gongyo, \textit{Log pluricanonical representations and abundance conjecture}, Compos. Math. \textbf{150} (2014), no. 4, 593--620.

%\bibitem[FH21]{FH21} O. Fujino and K. Hashizume, \textit{Existence of log canonical modifications and its applications}, arXiv: 2103.01417.

%\bibitem[FM00]{FM00} O. Fujino and S. Mori, \textit{A canonical bundle formula}, J. Differential Geom. \textbf{56} (2000), no. 1, 167--188.


%\bibitem[Fuk96]{Fuk96} S. Fukuda, \textit{On base point free theorem}, Kodai Math. J.19 (1996), no. 2,191--199.

%\bibitem[Gon11]{Gon11} Y. Gongyo, \textit{On the minimal model theory for dlt pairs of numerical kodaira dimension zero}, Math. Rest. Lett. \textbf{18} (2011), no. 5, 991--1000.

\bibitem[HH19]{HH19}  C. D. Hacon and J. Han, \textit{On a connectedness principle of Shokurov-Koll\'ar type,} Sci. China Math. \textbf{62} (2019), no. 3, 411--416.

\bibitem[HL21]{HL21} C. D. Hacon and J. Liu, \textit{Existence of flips for generalized lc pairs}, arXiv: 2105.13590v3.

%\bibitem[HMX14]{HMX14} C. D. Hacon, J. M\textsuperscript{c}Kernan, and C. Xu, \textit{ACC for log canonical thresholds}, Ann. of Math. \textbf{180} (2014), no. 2, 523--571.

%\bibitem[HM20]{HM20} C. D. Hacon and J. Moraga, \textit{On weak Zariski decompositions and termination of flips}, Math. Res. Lett. \textbf{27} (2020), no. 5, 1393--1421.


%\bibitem[HX13]{HX13} C. D. Hacon and C. Xu, \textit{Existence of log canonical closures}, Invent. Math. \textbf{192} (2013), no. 1, 161--195.

%\bibitem[HX15]{HX15} C. D. Hacon and C. Xu, \textit{Boundedness of log Calabi-Yau pairs of Fano type}, Math. Res. Lett. \textbf{22} (2015), no. 6, 1699--1716.

%\bibitem[HX16]{HX16} C. D. Hacon and C. Xu, \textit{On finiteness of B-representations and semi-log canonical abundance} in Minimal Models and Extremal Rays (Kyoto, 2011), Adv. Stud. Pure Math. \textbf{70} (2016), Math. Soc. Japan, Tokyo, 361--378. 

\bibitem[HL18]{HL18} J. Han and Z. Li, \textit{Weak Zariski decompositions and log terminal models for generalized polarized pairs}, arXiv: 1806.01234v2.

%\bibitem[HL20a]{HL20a} J. Han and Z. Li, \textit{On Fujita’s conjecture for pseudo-effective thresholds}, Math. Res. Lett. \textbf{27} (2020), no. 2, 377--396.

%\bibitem[HL20b]{HL20b} J. Han and Z. Li, \textit{On accumulation points of pseudo-effective thresholds}, Manuscripta math (2020).

%\bibitem[HL20c]{HL20c} J. Han and J. Liu, \textit{Effective birationality for sub-pairs with real coefficients}, arXiv: 2007.01849v1.

\bibitem[HLS19]{HLS19} J. Han, J. Liu, and V. V. Shokurov, \textit{ACC for minimal log discrepancies of exceptional singularities}, arXiv: 1903.04338v2.

\bibitem[HL19]{HL19} J. Han and W. Liu, \textit{On a generalized canonical bundle formula for generically finite morphisms}, arXiv: 1905.12542v3,  to appear in Ann. Inst. Fourier (Grenoble).

\bibitem[HL20]{HL20} J. Han and W. Liu, \textit{On numerical nonvanishing for generalized log canonical pairs}, Doc. Math. \textbf{25} (2020), 93--123.

%\bibitem[Has18]{Has18} K. Hashizume, \textit{Minimal model theory for relatively trivial log canonical pairs}, Ann. Inst. Fourier (Grenoble) \textbf{68} (2018), no. 5, 2069--2107.

%\bibitem[Has19]{Has19} K. Hashizume, \textit{Remarks on special kinds of the relative log minimal model program}, Manuscripta Math. \textbf{160} (2019), no. 3, 285--314.

%%\bibitem[Has20a]{Has20a} K. Hashizume, \textit{Finiteness of log abundant log canonical pairs in log minimal model program with scaling}, arXiv: 2005.12253v3.

\bibitem[Has20]{Has20} K. Hashizume, \textit{Non-vanishing theorem for generalized log canonical pairs with a polarization}, arXiv: 2012.15038v1.

\bibitem[HH20]{HH20}  K. Hashizume and Z. Hu, \textit{On minimal model theory for log abundant lc pairs}, J. Reine Angew. Math., \textbf{767} (2020), 109--159. 

\bibitem[Hu20]{Hu20} Z. Hu, \textit{Log abundance of the moduli b-divisors for lc-trivial fibrations}, arXiv: 2003.14379v3.

%\bibitem[Hu21]{Hu21} Z. Hu, \textit{An abundance theroem for generalised pairs}, arXiv: 2103.11813v1.

%\bibitem[Jia21]{Jia21} J. Jiao, \textit{On the Boundedness of Canonical Models}, arXiv: 2103.13609v1.

%\bibitem[Kaw84]{Kaw84} Y. Kawamata, \textit{The cone of curves of algebraic varieties}, Ann. of Math. \textbf{119} (1984), 603--633.

%\bibitem[Kaw92]{Kaw92} Y. Kawamata, \textit{Termination  of  log  flips  for  algebraic $3$-folds}, Internat. J. Math. \textbf{3} (1992), no. 5, 653--659.

\bibitem[Kaw98]{Kaw98} Y. Kawamata, \textit{Subadjunction of log canonical divisors, II}, Amer. J. Math. \textbf{120} (1998), 893--899.

%\bibitem[Kaw15]{Kaw15} Y. Kawamata, \textit{Variation of mixed Hodge structures and the positivity for algebraic fiber spaces}, Advanced Studies in Pure Mathematics, \textbf{65} (2015), 27--57.

%\bibitem[KMM87]{KMM87} Y. Kawamata, K. Matsuda, and K. Matsuki, \textit{Introduction to the minimal model problem}, Algebraic geometry, Sendai, 1985, 283--360, Adv. Stud. Pure Math., \textbf{10}, North-Holland, Amsterdam, 1987.

\bibitem[KMM94]{KMM94} S. Keel, K. Matsuki, and J. M\textsuperscript{c}Kernan, \textit{Log abundance theorem for threefolds}, Duke Math. J. \textbf{75} (1994), 99--119.



\bibitem[Kod63]{Kod63} K. Kodaira, \textit{On compact analytic surfaces}, II, III. Ann. of Math. (2) \textbf{77} (1963), 563--626, ibid., 78: 1--40.

%\bibitem[Kol84]{Kol84} J. Koll\'ar, \textit{The cone theorem}, Ann. of Math., \textbf{120} (1984), 1--5.

%\bibitem[Kol07]{Kol07} J. Koll\'ar, \textit{“Kodaira’s canonical bundle formula and adjunction}. In: \textit{Flips for 3-folds and 4-folds}. Ed. by A. Corti. Vol. 35. Oxford Lecture Series in Mathematics and its Applications. Oxford: Oxford University Press, 2007. Chap. 8, 134--162.

%\bibitem[Kol14]{Kol14} J. Koll\'ar, \textit{Semi-Normal Log Centres and Deformations of Pairs}, Proc. Edinburgh Math. Soc., \textbf{57} (2014), no. 1, 191--199.



\bibitem[KM98]{KM98} J. Koll\'{a}r and S. Mori, \textit{Birational geometry of algebraic varieties}, Cambridge Tracts in Math. \textbf{134} (1998), Cambridge Univ. Press.

\bibitem[LP20a]{LP20a} V. Laz\'ic and T. Peternell, \textit{On generalised abundance, I}, Publ. Res. Inst. Math. Sci. \textbf{56} (2020), no. 2, 353--389.

\bibitem[LP20b]{LP20b} V. Laz\'ic and T. Peternell, \textit{On generalised abundance, II}, Peking Mathematical Journal \textbf{3} (2020), 1--46.

%\bibitem[LMT20]{LMT20} V. Lazi\'c, J. Moraga, and N. Tsakanikas, \textit{Special termination for log canonical pairs}, arXiv: 2007.06458v1.

%\bibitem[LT19]{LT19} V. Lazi\'c and N. Tsakanikas, \textit{On the existence of minimal models for log canonical pairs}, to appear in Publ. Res. Inst. Math. Sci., arXiv: 1905.05576v3.

%\bibitem[LT21]{LT21} V. Lazi\'c and N. Tsakanikas, \textit{Special MMP for log canonical generalised pairs}, personal communication.

%\bibitem[Li20]{Li20} Z. Li, \textit{Boundedness of the base varieties of certain fibrations}, arXiv: 2002.06565v2.

%\bibitem[Li21]{Li21} Z. Li, \textit{Fujita’s conjecture on iterated accumulation points of pseudo-effective thresholds}, Selecta Mathematica \textbf{27} (2021), no. 9.

%\bibitem[Liu21]{Liu21} J. Liu, \textit{Sarkisov program for generalized pairs}, Osaka J. Math., \textbf{58} (2021), no. 4.

%\bibitem[LX21]{LX21} J. Liu and L. Xie, \textit{Number of singular points on projective surfaces}, arXiv: 2103.04522v1.

%\bibitem[Mor18]{Mor18} J. Moraga, \textit{Termination of pseudo-effective4-fold flips}, arXiv: 1802.10202v3.

%\bibitem[Nak86]{Nak86} N. Nakayama, \textit{Invariance of the plurigenera of algebraic varieties under minimal model conjectures}, Topology \textbf{25} (1986), no. 2, 237--251.

%\bibitem[Nak16]{Nak16} Y. Nakamura, \textit{On minimal log discrepancies on varieties with fixed Gorenstein index}. Michigan Math. J., 65 (1), 165--187, 2016.

\bibitem[Nak04]{Nak04} N. Nakayama, \textit{Zariski-decomposition and abundance}, MSJ Memoirs, vol. 14, Mathematical Society of Japan, Tokyo, 2004.

\bibitem[PS09]{PS09} Y. G. Prokhorov, V.V. Shokurov, \textit{Towards the second main theorem on complements}, J. Algebraic Geom., \textbf{18} (2009), 151--199.

%\bibitem[Sho96]{Sho96} V.V. Shokurov, \textit{3-fold log models}, J. Math. Sci. \textbf{81} (1996), no. 3, 2667--2699.

%\bibitem[Sho20]{Sho20} V.V. Shokurov, \textit{Existence and boundedness of $n$-complements}, arXiv: 2012.06495v1.

\end{thebibliography}
\end{document}